\documentclass[a4paper,oneside,10pt]{article}%
\usepackage{amsmath}
\usepackage{amsfonts}
\usepackage{amsthm}
\usepackage{indentfirst}
\usepackage{makeidx}
\usepackage{amsmath}
\usepackage{amssymb}
\usepackage{amsfonts}
\usepackage[perpage,symbol*]{footmisc}
\usepackage{bm}
\usepackage{graphicx}
\usepackage[T1]{fontenc}
\usepackage{titlesec}
\usepackage[tight,hang]{subfigure}
\usepackage{graphicx}%
\providecommand{\U}[1]{\protect\rule{.1in}{.1in}}
\providecommand{\U}[1]{\protect \rule{.1in}{.1in}}

\makeatletter

\newcommand{\Rmnum}[1]{\expandafter \@slowromancap \romannumeral #1@}
\makeatother

\newtheorem{theorem}{Theorem}[section]
\newtheorem{lemma}{Lemma}[section]

\newtheorem{definition}{Definition}[section]
\newtheorem{remark}{Remark}[section]

\numberwithin{equation}{section}

\begin{document}

\title{A maximum principle for controlled time-symmetric forward-backward doubly
stochastic differential equation with initial-terminal sate
constraints\thanks{This work was supported by National Natural Science
Foundation of China (No. 11171187, No. 10871118 and No. 10921101); supported
by the Programme of Introducing Talents of Discipline to Universities of China
(No. B12023); supported by Program for New Century Excellent Talents in
University of China.}}
\author{Shaolin Ji\thanks{Institute for Financial Studies and Institute of
Mathematics, Shandong University, Jinan, Shandong 250100, PR China
(jsl@sdu.edu.cn, Fax: +86 0531 88564100).}
\and Qingmeng Wei\thanks{Institute of mathematics, Shandong University, Jinan,
Shandong 250100, PR China. (qingmengwei@gmail.com)}
\and Xiumin Zhang\thanks{Institute of mathematics, Shandong University, Jinan,
Shandong 250100, PR China. }}
\maketitle

\begin{abstract}
In this paper, we study the optimal control problem of a controlled
time-symmetric forward-backward doubly stochastic differential equation with
initial-terminal sate constraints. Applying the terminal perturbation method
and Ekeland's variation principle, a necessary condition of the stochastic
optimal control, i.e., stochastic maximum principle is derived. Applications
to backward doubly stochastic linear-quadratic control models are investigated.
\end{abstract}

\textbf{Keywords}
Time-symmetric forward-backward doubly stochastic differential
equations, Ekeland's variation principle,
State constraints, Stochastic maximum
principle.

\textbf{AMS}
93E20, 60H10

\section{Introduction}

It is well known that a general coupled forward-backward stochastic
differential equations (FBSDEs for short) consists of a forward SDE of
It\^{o}'s type and a backward SDE of Pardoux-Peng's (see
\cite{El.Karoui-peng97, Pardoux-Peng90}). Since Antonelli \cite{Antonelli}
first studied FBSDEs in early 1990s, FBSDEs have been studied widely in many
papers (see \cite{Hu-Peng, Ma-Protter-Yong, Ma-Yong, Yong}). FBSDEs are often
encountered in the optimization problem when applying stochastic maximum
principle (see \cite{Hu1991, Yong-Zhou}). In finance, FBSDEs are used when
considering problems with the large investors, see \cite{Buckdahn-Hu1998,
Cvitanic-Ma, Ma-Yong}. Such equations are also used in the potential theory
(see \cite{Hu1999}). Moreover, one can apply FBSDEs to study Homogenization
and singular perturbation of certain quasilinear parabolic PDEs with periodic
structures (see \cite{Briand-Hu1999, Buckdahn-Hu1998-1}).

In order to produce a probabilistic representation of certain quasilinear
stochastic partial differential equations (SPDEs for short), Pardoux and Peng
\cite{Pardoux-Peng98} first introduced backward doubly stochastic differential
equations (BDSDEs for short) and proved the existence and uniqueness theorem
of BDSDEs. Using such BDSDEs they proved the existence and uniqueness theorem
of those quasilinear SPDEs and thus significantly extended the famous
Feynman-Kac formula for such SPDEs.

Peng and Shi \cite{Peng-Shi03} studied the following time-symmetric
forward-backward doubly stochastic differential equations (FBDSDEs for short):%
\[
\left\{
\begin{array}
[c]{rrl}%
-dx_{t} & = & F(t,x_{t},z_{t},y_{t},q_{t},u_{t})dt+G(t,x_{t},z_{t},y_{t}%
,q_{t},u_{t})dW_{t}-z_{t}dB_{t},\ 0\leq t\leq T,\\
x_{0} & = & \xi,\\
-dy_{t} & = & f(t,x_{t},z_{t},y_{t},q_{t},u_{t})dt+g(t,x_{t},z_{t},y_{t}%
,q_{t},u_{t})dB_{t}-q_{t}dW_{t},\ 0\leq t\leq T,\\
y_{T} & = & \eta,
\end{array}
\right.
\]
which generalized the general FBSDEs. Here the forward equation is ``forward"
with respect to a standard stochastic integral $dW_{t}$, as well as
``backward" with respect to a backward stochastic integral $dB_{t}$; the
coupled ``backward equation" is ``forward" under the backward stochastic
integral $dB_{t}$ and ``backward" under the forward one. In other wards, both
the forward equation and the backward one are BDSDEs with different directions
of stochastic integral. Under certain monotonicity conditions, they proved the
uniqueness and existence theorem for these equations. In \cite{Han-Peng-Wu},
when deriving the stochastic maximum principle of backward doubly stochastic
optimal control problems, Han, Peng and Wu showed that this kind equations are
just the state equation and adjoint equation of their optimal control problem.

In this paper, we study a stochastic optimal control problem with
initial-terminal state constraints where the controlled system is described by
the above time-symmetric FBDSDEs. We suppose that the initial state $\xi$ and
the terminal state $\eta$ fall in two convex sets, respectively, and the
corresponding states $x_{T}^{(\xi,\eta,u(\cdot))}$ and $y_{0}^{(\xi
,\eta,u(\cdot))}$ satisfy the constraints $E(\psi(x_{T}^{(\xi,\eta,u(\cdot
))}))=a$ and $E(h(y_{0}^{(\xi,\eta,u(\cdot))}))=b$ respectively. Then we
minimize the following cost function:%
\[
J(\xi,\eta,u(\cdot))\triangleq E[\int_{0}^{T}%
l(x(t),z(t),y(t),q(t),u(t),t)dt+\chi(\xi)+\lambda(\eta)+\phi(x(T))+\gamma
(y(0))].
\]

It is well-known that the maximum principle is an important approach to study
optimal control problems. The systematic account on this theory can be found
in \cite{Bensoussan, Yong-Zhou}. When the controlled system under
consideration is assumed to be with state constraints, especially with
sample-wise constraints, the corresponding stochastic optimal control problems
are difficult to solve. A sample-wise constraint requires that the state be in
a given set with probability $1$; for example, a nonnegativity constraint on the
wealth process, i.e., bankruptcy prohibition in financial markets. In order to
deal with such optimal control problems, an approach named ``terminal
perturbation method" was introduced and applied in financial optimization
problems recently (see \cite{Ji, Ji-Peng, Ji-Zhou, Ji-Zhou10}). This method is
based on the dual method or martingale method introduced by Bieleckiet in
\cite{Bielecki-Pliska-Zhou} and El Karoui, Peng and Quenez in
\cite{EL.Karoui-peng01}. It mainly applies Ekeland's variational principle to
tackle the state constraints and derive a stochastic maximum principle which
characterizes the optimal solution. For other works about the optimization
problem with state constraints, the readers may refer to \cite{Yong10-1,
Yong10-2}. In this paper, a stochastic maximum principle is obtained for the
controlled time-symmetric FBDSDEs with initial-terminal state constraints by
using Ekeland's variational principle.

We give three specific applications to illustrate our theoretical results. In
the first application, the controlled state equations are composed of a normal
FSDE and a BDSDE. By introducing a backward formulation of the controlled
system (inspired by \cite{Ji-Zhou}), we present the stochastic maximum
principle for the optimal control. As a special case, we only consider one
BDSDE as our state equation in the second application. As stated in the last
application, our results can be applied in forward-backward doubly stochastic
linear-quadratic (LQ) optimal control problems. The explicit expression of the
optimal control is derived. Since the control system of SPDEs can be
transformed to the relevant control system of FBDSDEs, our results can be used
to solve the optimal control problem of one kind of SPDEs.

This paper is organized as follows. In section 2.1, we recall some
preliminaries. And we formulate our control problem in section 2.2. In seciton
2.3, by applying Ekeland's variation principle we obtain a stochastic maximum
principle of this controlled time-symmetric FBDSDEs with initial-terminal
state constraints. Some applications are given in the last section.

\section{The main problem}

\subsection{Preliminaries}

Let us first recall the existence and uniqueness results of the BDSDE which
was introduced by Pardoux and Peng \cite{Pardoux-Peng98}, and an extension of
the well-known It\^{o}'s formula which would be often used in this paper.

Let $(\Omega,%
\mathcal{F}%
,P)$ be a probability space, and $T>0$ be fixed throughout this paper. Let
$\{\mathbf{W}_{t},0\leq t\leq T\}$ and $\{\mathbf{B}_{t},0\leq t\leq T\}$ be
two mutually independent standard Brownian motion processes, with values in
$\mathbb{R}^{d},\ \mathbb{R}^{l}$, respectively, defined on $(\Omega,
\mathcal{F}
,P)$. Let $
\mathbb{N}
$ denote the class of $P$-null set of $
\mathcal{F}
$. For each $t\in\lbrack0,T]$, we define: $
\mathcal{F}
_{t}\triangleq
\mathcal{F}
_{t}^{W}\vee
\mathcal{F}
_{t,T}^{B}$, where
\[
\mathcal{F}_{t}^{W}  = 
\mathcal{F}_{0,t}^{W}  = \sigma\{W_{r}-W_{0};0\leq r\leq t\}\vee\mathbb{N},\
\mathcal{F}_{t}^{B}  = 
\mathcal{F}_{0,t}^{B}  = \sigma\{B_{r}-B_{t};t\leq r\leq T\}\vee\mathbb{N}.
\]
Note that the collection $\{
\mathcal{F}_{t},t\in\lbrack0,T]\}$ is neither increasing nor decreasing, and it does not
constitute a filtration.

For any Euclidean space $H$, we denote by $\langle\cdot,\cdot\rangle$ the
scale product of $H$. The Euclidean norm of a vector $y\in%
\mathbb{R}
^{k}$ will be denoted by $|y|,$ and for a $d\times n$ matrix A, we define
$||A||=\sqrt{Tr(AA^{\ast})}.$

For any $n\in N$, let $M^{2}(0,T;%
\mathbb{R}
^{n})$ denote the set of (classes of $dP\otimes dt$ a.e. equal) $n$%
-dimensional jointly measurable stochastic processes $\{\varphi_{t}%
;t\in\lbrack0,T]\}$ which satisfy:

(i) $E\int_{0}^{T}|\varphi_{t}|^{2}dt<\infty;$ 
(ii) $\varphi_{t}$ is $\mathcal{F}_{t}$-measurable, for a.e. $t\in\lbrack0,T].$

We denote by $S^{2}(0,T;\mathbb{R}^{n})$ the set of continuous $n$-dimensional stochastic processes which satisfy:

(i) $E(\sup\limits_{0\leq t\leq T}|\varphi_{t}|^{2})<\infty;$ 
(ii) $\varphi_{t}$ is $%
\mathcal{F}%
_{t}$-measurable, for any $t\in\lbrack0,T].$

Let
\[%
\begin{array}
[c]{l}%
f:\Omega\times\lbrack0,T]\times%
\mathbb{R}
^{k}\times%
\mathbb{R}
^{k\times d}\rightarrow%
\mathbb{R}
^{k},\
g:\Omega\times\lbrack0,T]\times%
\mathbb{R}
^{k}\times%
\mathbb{R}
^{k\times d}\rightarrow%
\mathbb{R}
^{k\times l},
\end{array}
\]
be jointly measurable and such that for any $(y,q)\in\mathbb{R}^{k}%
\times\mathbb{R}^{k\times d},$
$
f(\cdot,y,q) \in M^{2}(0,T;\mathbb{R}^{k}),\ g(\cdot,y,q) \in M^{2}%
(0,T;\mathbb{R}^{k\times l}).
$

Moreover, we assume that there exist constants $C>0$ and $0<\alpha<1$ such
that for any $(\omega,t)\in\Omega\times\lbrack0,T],\ (y_{1},q_{1}),(y_{2}
,q_{2})\in\mathbb{R}^{k}\times\mathbb{R}^{k\times l},$
\begin{equation}
\begin{array}
[c]{lll}
&  & |f(t,y_{1},q_{1})-f(t,y_{2},q_{2})|^{2} \leq C(|y_{1}-y_{2}|^{2}%
+\|q_{1}-q_{2}\|^{2});\\
&  & \|g(t,y_{1},q_{1})-g(t,y_{2},q_{2})\|^{2} \leq C|y_{1}-y_{2}|^{2}%
+\alpha\|q_{1}-q_{2}\|^{2}.
\end{array}
\tag{H}%
\end{equation}

Given $\eta\in L^{2}(\Omega, \mathcal{F} _{T},P;\mathbb{R}^{k})$, we consider
the following BDSDE:
\begin{equation}
y_{t}=\eta+\int_{t}^{T}f(s,y_{s},q_{s})ds+\int_{t}^{T}g(s,y_{s},q_{s}%
)dB_{s}-\int_{t}^{T}q_{s}dW_{s},\ 0\leq t\leq T. \tag{2.1}%
\end{equation}

We note that the integral with respect to $\{B_{t}\}$ is a ``backward It\^{o}
integral" and the integral with respect to $\{W_{t}\}$ is a standard forward
It\^{o} integral. These two types of integrals are particular cases of the
It\^{o}-Skorohod integral, see Nualart and Pardoux \cite{Nualart-Pardoux}.

By Theorem 1.1 in \cite{Pardoux-Peng98}, the above equation (2.1) has a unique
solution$\ (y,q)\in S^{2}(0,T;%
\mathbb{R}
^{k})\times M^{2}(0,T;%
\mathbb{R}
^{k\times d})$.

Next let us recall an extension of the well-known It\^{o}'s formula in [17]
which would be often used in this paper.

\begin{lemma}
Let $\alpha\in S^{2}(0,T; \mathbb{R}^{k}),$ $\beta\in M^{2}(0,T;\mathbb{R}%
^{k}),$ $\gamma\in M^{2}(0,T;\mathbb{R}^{k\times l}),$ $\delta\in
M^{2}(0,T;\mathbb{R}^{k\times d})$ be such that:%
\[
\alpha_{t}=\alpha_{0}+\int_{0}^{t}\beta_{s}ds+\int_{0}^{t}\gamma_{s}%
dB_{s}+\int_{0}^{t}\delta_{s}dW_{s},\ 0\leq t\leq T.
\]
Then,
\[
\begin{array}
[c]{rrl}%
|\alpha_{t}|^{2} & = & |\alpha_{0}|^{2}+2\int_{0}^{t}(\alpha_{s},\beta
_{s})ds+2\int_{0}^{t}(\alpha_{s},\gamma_{s}dB_{s})-\int_{0}^{t}||\gamma
_{s}||^{2}ds+\int_{0}^{t}||\delta_{s}||^{2}ds+2\int_{0}^{t}(\alpha_{s},\delta_{s}dW_{s}),\\
E|\alpha_{t}|^{2} & = & E|\alpha_{0}|^{2}+2E\int_{0}^{t}(\alpha_{s},\beta
_{s})ds-E\int_{0}^{t}||\gamma_{s}||^{2}ds+E\int_{0}^{t}||\delta_{s}||^{2}ds.
\end{array}
\]
Generally, for $\phi\in C^{2}(\mathbb{R}^{k}),$
\[
\begin{array}
[c]{rl}
\phi(\alpha_{t})= & \phi(\alpha_{0})+\int_{0}^{t}(\phi^{\prime}(\alpha
_{s}),\beta_{s})ds+\int_{0}^{t}(\phi^{\prime}(\alpha_{s}),\gamma_{s}%
dB_{s})+\int_{0}^{t}(\phi^{\prime}(\alpha_{s}),\delta_{s}dW_{s})\\
& -\frac{1}{2}\int_{0}^{t}Tr[\phi^{\prime\prime}(\alpha_{s})\gamma_{s}%
\gamma_{s}^{\ast}]ds+\frac{1}{2}\int_{0}^{t}Tr[\phi^{\prime\prime}(\alpha
_{s})\delta_{s}\delta_{s}^{\ast}]ds.
\end{array}
\]
\end{lemma}

\subsection{Problem formulation}

Let\ $K$ be a nonempty convex subset of $%
\mathbb{R}
^{n\times d}$. We set{\LARGE \ }%
\[
U[0,T]=\{u(\cdot)|u(t)\in K,\ a.e.,a.s.,0\leq t\leq T;u(\cdot)\in M^{2}(0,T;%
\mathbb{R}
^{n\times d})\}.
\]
An element of $U[0,T]$ is called an admissible control. Now let%
\[%
\begin{array}
[c]{l}%
F:\Omega\times\lbrack0,T]\times
\mathbb{R}
^{n}\times
\mathbb{R}
^{n\times l}\times
\mathbb{R}
^{k}\times
\mathbb{R}
^{k\times d}\times
\mathbb{R}
^{n\times d}\rightarrow
\mathbb{R}
^{n},\\
G:\Omega\times\lbrack0,T]\times
\mathbb{R}
^{n}\times
\mathbb{R}
^{n\times l}\times
\mathbb{R}
^{k}\times
\mathbb{R}
^{k\times d}\times
\mathbb{R}
^{n\times d}\rightarrow
\mathbb{R}
^{n\times d},\\
f:\Omega\times\lbrack0,T]\times\mathbb{R}^{n}\times\mathbb{R}^{n\times l}\times\mathbb{R}^{k}\times\mathbb{R}^{k\times d}\times\mathbb{R}^{n\times d}\rightarrow\mathbb{R}^{k},\\
g:\Omega\times\lbrack0,T]\times\mathbb{R}^{n}\times\mathbb{R}^{n\times l}\times\mathbb{R}^{k}\times\mathbb{R}^{k\times d}\times\mathbb{R}^{n\times d}\rightarrow\mathbb{R}^{k\times l},
\end{array}
\]
be jointly measurable such that for any $(x,z,y,q)\in%
\mathbb{R}
^{n}\times%
\mathbb{R}
^{n\times l}\times%
\mathbb{R}
^{k}\times%
\mathbb{R}
^{k\times d}$ and any $u(\cdot)$ $\in U[0,T]$%
\[%
\begin{array}
[c]{lll}%
F(\cdot,x,z,y,q,u(\cdot)) \in M^{2}(0,T; \mathbb{R}^{n}), \ G(\cdot
,x,z,y,q,u(\cdot)) \in M^{2}(0,T;\mathbb{R} ^{n\times d}), &  & \\
f(\cdot,x,z,y,q,u(\cdot)) \in M^{2}(0,T;\mathbb{R}^{k}), \ g(\cdot
,x,z,y,q,u(\cdot)) \in M^{2}(0,T;\mathbb{R}^{k\times l}). &  &
\end{array}
\]
Let
\[
\zeta(t)  =  (x(t),z(t),y(t),q(t))^{T},\
A(t,\zeta) =  (-F,-G,-f,-g)^{T}(t,\zeta).
\]
We assume

(H1) $\forall\zeta_{1}=(x_{1},z_{1},y_{1},q_{1})$, $\zeta_{2}=(x_{2}%
,z_{2},y_{2},q_{2})\in%
\mathbb{R}
^{n}\times%
\mathbb{R}
^{n\times l}\times%
\mathbb{R}
^{k}\times%
\mathbb{R}
^{k\times d}$ and $t\in\lbrack0,T]$, there exists a constant $\mu>0$ such that
the following monotonicity condition holds for any $u(\cdot)$ $\in U[0,T]$%
\[
\langle A(t,\zeta_{1})-A(t,\zeta_{2}),\zeta_{1}-\zeta_{2}\rangle\leq-\mu
|\zeta_{1}-\zeta_{2}|^{2}.
\]

(H2) There exist constants $C>0$ and $0<\alpha<\frac{1}{2}$ such that for any
$(\omega,t)\in\Omega\times\lbrack0,T],$ $u\in\mathbb{R}^{n\times d}$,
$(x_{1},z_{1},y_{1},q_{1}),(x_{2},z_{2},y_{2},q_{2})\in\mathbb{R}^{n}%
\times\mathbb{R}^{n\times l}\times\mathbb{R}^{k}\times\mathbb{R}^{k\times d}$
the following conditions hold:
\[%
\begin{array}
[c]{rl}
& |\vartheta(t,x_{1},z_{1},y_{1},q_{1},u)-\vartheta(t,x_{2},z_{2} ,y_{1}%
,q_{1},u)|^{2}
\leq  C(|x_{1}-x_{2}|^{2}+\|z_{1}-z_{2}\|^{2}+|y_{1}-y_{2}|^{2}+\|q_{1}
-q_{2}\|^{2}),\\
& \|\varkappa(t,x_{1},z_{1},y_{1},q_{1},u)-\varkappa(t,x_{2},z_{2}
,y_{1},q_{1},u)\|^{2}
\leq  C(|x_{1}-x_{2}|^{2}+|y_{1}-y_{2}|^{2})+\alpha(\|z_{1}-z_{2}%
\|^{2}+\|q_{1}-q_{2}\|^{2}),
\end{array}
\]
where $\vartheta=(F,f),$ $\varkappa=(G,g).$

(H3) $F,G,f,g,\psi,h,l,\chi,\lambda,\phi$ and $\gamma$ are continuous in their
arguments and continuously differentiable in $(x,z,y,q,u)$, and the
derivatives of $F,G,f,g$ in $(x,z,y,q,u)$ are bounded and $0<\|g_{z}%
(\cdot)\|<\frac{1}{2},$ $0<\|G_{z}(\cdot)\|<\frac{1}{2},$ $0<\|g_{q}%
(\cdot)\|<\frac{1}{2},$ $0<\|G_{q}(\cdot)\|<\frac{1}{2};$ the derivatives of
$l$ in $(x,y,z,q,u)$ are bounded by $C(1+|x|+|z|+|y|+|q|+\|u\|)$, and the
derivatives of $\phi$, $\chi$ and $\psi$ in $x$ are bounded by $C(1+|x|)$ ;
$\gamma$, $\lambda$ and $h$ in $y$ are bounded by $C(1+|y|)$.

Given $\xi\in L^{2}(\Omega,\mathcal{F}_{0},P;\mathbb{R}^{n})$, $\eta\in
L^{2}(\Omega,\mathcal{F}_{T},P;\mathbb{R}^{k})$ and $\forall u(\cdot)\in
U[0,T],$ let us consider the following time-symmetric FBDSDE:
\begin{equation}
\left\{
\begin{array}
[c]{rrl}%
-dx_{t} & = & F(t,x_{t},z_{t},y_{t},q_{t},u_{t})dt+G(t,x_{t},z_{t},y_{t}%
,q_{t},u_{t})dW_{t}-z_{t}dB_{t},\\
x_{0} & = & \xi,\\
-dy_{t} & = & f(t,x_{t},z_{t},y_{t},q_{t},u_{t})dt+g(t,x_{t},z_{t},y_{t}%
,q_{t},u_{t})dB_{t}-q_{t}dW_{t},\\
y_{T} & = & \eta.
\end{array}
\right.  \tag{2.2}%
\end{equation}
Recall Theorem 2.2 in \cite{Peng-Shi03}. We have

\begin{theorem}
For given $\xi\in L^{2}(\Omega,\mathcal{F}_{0},P;\mathbb{R}^{n})$, $\eta\in
L^{2}(\Omega,\mathcal{F}_{T},P;\mathbb{R}^{k})$ and $\forall u(\cdot)\in
U[0,T]$, assume (H1)$\sim$(H3), then (2.2) exists a unique $\mathcal{F}_{t}%
$-adapted solution $(x(t),z(t),$ $y(t),q(t))$.
\end{theorem}

In (2.2), we regard $\xi,$ $\eta,$ $u(\cdot)$ as controls. $\xi,$ $\eta,$
$u(\cdot)$ can be chosen from the following admissible set :%
\[
\begin{array}
[c]{rl}%
U= & \{(\xi,\eta,u(\cdot))|\xi\in K_{1}\subset
\mathbb{R}
^{n}a.s.\text{, }\eta\in K_{2}\subset
\mathbb{R}
^{k}\ a.s.,\ 
E[|\xi|^{2}]<\infty,E[|\eta|^{2}]<\infty\text{, }u(\cdot)\in U[0,T]\},
\end{array}
\]
where $K_{1}$ and $K_{2}$ are convex.

\begin{remark}
A main assumption in this paper is the control domains are convex. For the
terminal perturbation method, it is difficult to weaken or completely remove
these assumptions. Until now, it remains an interesting and challenging open problem.
\end{remark}

We also assume the state constraints
\[
E(\psi(x_{T}^{(\xi,\eta,u(\cdot))}))=a,\ E(h(y_{0}^{(\xi,\eta,u(\cdot))}))=b.
\]
For each $(\xi,\eta,u(\cdot))\in U$, consider the following cost function:%
\begin{equation}%
\begin{array}
[c]{lll}%
J(\xi,\eta,u(\cdot)) \triangleq E[\int_{0}^{T}%
l(x(t),z(t),y(t),q(t),u(t),t)dt+\chi(\xi)+\lambda(\eta)+\phi(x(T)) +\gamma(y(0))]. &  &
\end{array}
\tag{2.3}%
\end{equation}
Our optimization problem is:%
\begin{equation}%
\begin{array}
[c]{cc}
& \inf\limits_{(\xi,\eta,u(\cdot))\in U}J(\xi,\eta,u(\cdot))\\
\text{subject to} & E(\psi(x_{T}^{(\xi,\eta,u(\cdot))}))=a\text{; }%
E(h(y_{0}^{(\xi,\eta,u(\cdot))}))=b.
\end{array}
\tag{2.4}%
\end{equation}

\begin{definition}
A triple of random variable $(\xi,\eta,u(\cdot))\in U$ is called feasible for
given $a\in\mathbb{R}^{n}$, $b\in\mathbb{R}^{k}$ if the solution (2.2) satisfy
$E(\psi(x_{T}^{(\xi,\eta,u(\cdot))}))=a$ and $E(h(y_{0}^{(\xi,\eta,u(\cdot
))}))=b$. We shall denote by $\mathbb{N}(a,b)$ the set of all feasible
$(\xi,\eta,u(\cdot))$ for any given $a$ and $b$.
\end{definition}

A feasible $(\xi^{\ast},\eta^{\ast},u^{\ast}(\cdot))$ is called optimal if it
attains the minimum of $J(\xi,\eta,u(\cdot))$ over $%
\mathbb{N}
(a,b)$.

The aim of this paper is to obtain a characterization of $(\xi^{\ast}%
,\eta^{\ast},u^{\ast}(\cdot))$, i.e., the stochastic maximum principle.

\subsection{Stochastic Maximum Principle}

Using Ekeland's variational principle, we derive maximum principle for the
optimization problem (2.4) in this section. For simplicity, we first study the
case where $l(y(t),z(t),y(t),q(t),u(t),t)=0$, $\chi(x)=0$ and $\lambda(y)=0$
in subsection 2.3.1-2.3.3, and then present the results for the general case
in subsection 2.3.4.

\subsubsection{Variational equations}

For $(\xi^{1},\eta^{1},u^{1}(\cdot)),(\xi^{2},\eta^{2},u^{2}(\cdot))\in U,$ we
define a metric in $U$ by%
\[
d(\xi^{1},\eta^{1},u^{1}(\cdot)),(\xi^{2},\eta^{2},u^{2}(\cdot))\triangleq
(E|\xi^{1}-\xi^{2}|^{2})^{\frac{1}{2}}+E|\eta^{1}-\eta^{2}|^{2})^{\frac{1}{2}%
}+(\|u^{1}(\cdot)-u^{2}(\cdot)\|^{2})^{\frac{1}{2}}.
\]
It is obvious that $(U,d(\cdot,\cdot))$ is a complete metric space.

Let $(\xi^{\ast},\eta^{\ast},u^{\ast}(\cdot))$ be optimal and $(x^{\ast}%
(\cdot),z^{\ast}(\cdot),y^{\ast}(\cdot),q^{\ast}(\cdot)) $\ be the
corresponding state processes of (2.2). $\forall(\xi,\eta,u(\cdot))\in U$ and
$\forall0\leq\rho\leq1,$%
\[
(\xi^{\rho},\eta^{\rho},u^{\rho}(\cdot))\triangleq(\xi^{\ast}+\rho(\xi
-\xi^{\ast}),\eta^{\ast}+\rho(\eta-\eta^{\ast}),u^{\ast}(\cdot)+\rho
(u(\cdot)-u^{\ast}(\cdot)))\in U.
\]
Let $(x^{\rho}(\cdot),z^{\ \rho}(\cdot),y^{\rho}(\cdot),q^{\ \rho}(\cdot))$ be
the state processes of (2.2) associated with $(\xi^{\rho},\eta^{\rho},$
$u^{\rho}(\cdot))$.
\\
To derive the first-order necessary condition, we let $(\hat{x}(\cdot),\hat
{z}(\cdot),\hat{y}(\cdot),\hat{q}(\cdot))$ be the solution of the following
time-symmetric FBDSDE:%
\begin{equation}
\left\{
\begin{array}
[c]{rrl}%
-d\hat{x}(t) & = & (F_{x}^{\ast}(t)\hat{x}(t)+F_{z}^{\ast}(t)\hat{z}%
(t)+F_{y}^{\ast}(t)\hat{y}(t)+F_{q}^{\ast}(t)\hat{q}(t)+F_{u}^{\ast
}(t)(u(\cdot)-u^{\ast}(\cdot)))dt\\
&  & +(G_{x}^{\ast}(t)\hat{x}(t)+G_{z}^{\ast}(t)\hat
{z}(t)+G_{y}^{\ast}(t)\hat{y}(t)+G_{q}^{\ast}(t)\hat{q}(t) +G_{u}^{\ast}(t)(u(\cdot)-u^{\ast}(\cdot)))dW_{t}-\hat{z}(t)dB(t),\\
\hat{x}(0) & = & \xi-\xi^{\ast},\\
-d\hat{y}(t) & = & (f_{x}^{\ast}(t)\hat{x}(t)+f_{z}^{\ast}(t)\hat{z}%
(t)+f_{y}^{\ast}(t)\hat{y}(t)+f_{q}^{\ast}(t)\hat{q}(t)+f_{u}^{\ast
}(t)(u(\cdot)-u^{\ast}(\cdot)))dt\\
&  &+(g_{x}^{\ast}(t)\hat{x}(t)+g_{z}^{\ast}(t)\hat
{z}(t)+g_{y}^{\ast}(t)\hat{y}(t)+g_{q}^{\ast}(t)\hat{q}(t) +g_{u}^{\ast}(t)(u(\cdot)-u^{\ast}(\cdot)))dB_{t}-\hat{q}(t)dW(t),\\
\hat{y}(T) & = & \eta-\eta^{\ast},
\end{array}
\right.  \tag{2.5}%
\end{equation}
where$\ H_{k}^{\ast}(t)=H_{k}(t,x^{\ast}(\cdot),z^{\ast}(\cdot),y^{\ast}
(\cdot),q^{\ast}(\cdot),u^{\ast}(\cdot))$ for $H=F,G,f,g$ , $k=x,z,y,q,u,$
respectively. Equation (2.5) is called the variation equation.

Set
\begin{align}
\tilde{x}_{\rho}(t)    =\rho^{-1}[x_{\rho}(t)-x^{\ast}(t)]-\hat
{x}(t),\ 
\tilde{z}_{\rho}(t)   =\rho^{-1}[z_{\rho}(t)-z^{\ast}(t)]-\hat
{z}(t),\nonumber\\
\tilde{y}_{\rho}(t)    =\rho^{-1}[y_{\rho}(t)-y^{\ast}(t)]-\hat
{y}(t),\ \tag{2.6}
\tilde{q}_{\rho}(t)   =\rho^{-1}[q_{\rho}(t)-q^{\ast}(t)]-\hat
{q}(t).\nonumber
\end{align}
We have the following convergence.

\begin{lemma}
Assuming $(H1)\sim(H3)$ we have%
\begin{equation}%
\begin{array}
[c]{c}%
\lim\limits_{\rho\rightarrow0}\sup\limits_{0\leq t\leq T}E[|\tilde{x}_{\rho
}(t)|^{2}]=0,\
\lim\limits_{\rho\rightarrow0}E[\int_{0}^{T}\|\tilde{z}_{\rho}(t)\|^{2}%
dt]=0,\\
\lim\limits_{\rho\rightarrow0}\sup\limits_{0\leq t\leq T}E[|\tilde{y}_{\rho
}(t)|^{2}]=0,\
\lim\limits_{\rho\rightarrow0}E[\int_{0}^{T}\|\tilde{q}_{\rho}(t)\|^{2}dt]=0.
\end{array}
\tag{2.7}%
\end{equation}

\end{lemma}

\begin{proof}
From (2.2) and \ (2.5), we have%
\begin{equation}
\left\{
\begin{array}
[c]{lll}%
-d\tilde{y}_{\rho}(t) = \rho^{-1}[f(x_{\rho}(t),z_{\rho}(t),y_{\rho
}(t),q_{\rho}(t),u_{\rho}(t),t)-f(x^{\ast}(t),z^{\ast}(t),y^{\ast}(t),q^{\ast
}(t),u^{\ast}(t),t) &  & \\
\qquad\qquad-\rho f_{x}^{\ast}(t)\hat{x}(t)-\rho f_{z}^{\ast}(t)\hat
{z}(t)-\rho f_{y}^{\ast}(t)\hat{y}(t)-\rho f_{q}^{\ast}(t)\hat{q}(t)-\rho
f_{u}^{\ast}(t)(u(\cdot)-u^{\ast}(\cdot))]dt &  & \\
\qquad\qquad+\rho^{-1}[g(x_{\rho}(t),z_{\rho}(t),y_{\rho}(t),q_{\rho
}(t),u_{\rho}(t),t)-g(x^{\ast}(t),z^{\ast}(t),y^{\ast}(t),q^{\ast}(t),u^{\ast
}(t),t) &  & \\
\qquad\qquad-\rho g_{x}^{\ast}(t)\hat{x}(t)-\rho g_{z}^{\ast}(t)\hat
{z}(t)-\rho g_{y}^{\ast}(t)\hat{y}(t)-\rho g_{q}^{\ast}(t)\hat{q}(t)-\rho
g_{u}^{\ast}(t)(u(\cdot)-u^{\ast}(\cdot))]dB_{t} -\tilde{q}_{\rho}(t)dW_{t}, &  & \\
\tilde{y}_{\rho}(T) = 0. &  &
\end{array}
\right.  \tag{2.8}%
\end{equation}
Let%
\[%
\begin{array}
[c]{rrl}%
A^{\rho}(t) & = & \int_{0}^{1}f_{x}(x^{\ast}(t)+\lambda\rho(\hat{x}%
(t)+\tilde{x}_{\rho}(t)),z^{\ast}(t)+\lambda\rho(\hat{z}(t)+\tilde{z}_{\rho
}(t)),y^{\ast}(t)+\lambda\rho(\hat{y}(t)+\tilde{y}_{\rho}(t)),\\
&  & q^{\ast}(t)+\lambda\rho(\hat{q}(t)+\tilde{q}_{\rho}(t)),u^{\ast
}(t)+\lambda\rho(u(t)-u^{\ast}(t)),t)d\lambda,\\
B^{\rho}(t) & = & \int_{0}^{1}f_{z}(x^{\ast}(t)+\lambda\rho(\hat{x}%
(t)+\tilde{x}_{\rho}(t)),z^{\ast}(t)+\lambda\rho(\hat{z}(t)+\tilde{z}_{\rho
}(t)),y^{\ast}(t)+\lambda\rho(\hat{y}(t)+\tilde{y}_{\rho}(t)),\\
&  & q^{\ast}(t)+\lambda\rho(\hat{q}(t)+\tilde{q}_{\rho}(t)),u^{\ast
}(t)+\lambda\rho(u(t)-u^{\ast}(t)),t)d\lambda,\\
C^{\rho}(t) & = & \int_{0}^{1}f_{y}(x^{\ast}(t)+\lambda\rho(\hat{x}%
(t)+\tilde{x}_{\rho}(t)),z^{\ast}(t)+\lambda\rho(\hat{z}(t)+\tilde{z}_{\rho
}(t)),y^{\ast}(t)+\lambda\rho(\hat{y}(t)+\tilde{y}_{\rho}(t)),\\
&  & q^{\ast}(t)+\lambda\rho(\hat{q}(t)+\tilde{q}_{\rho}(t)),u^{\ast
}(t)+\lambda\rho(u(t)-u^{\ast}(t)),t)d\lambda,\\
D^{\rho}(t) & = & \int_{0}^{1}f_{q}(x^{\ast}(t)+\lambda\rho(\hat{x}%
(t)+\tilde{x}_{\rho}(t)),z^{\ast}(t)+\lambda\rho(\hat{z}(t)+\tilde{z}_{\rho
}(t)),y^{\ast}(t)+\lambda\rho(\hat{y}(t)+\tilde{y}_{\rho}(t)),\\
&  & q^{\ast}(t)+\lambda\rho(\hat{q}(t)+\tilde{q}_{\rho}(t)),u^{\ast
}(t)+\lambda\rho(u(t)-u^{\ast}(t)),t)d\lambda,\\
E^{\rho}(t) & = & [A^{\rho}(t)-f_{x}^{\ast}(t)]\hat{x}(t)+[B^{\rho}%
(t)-f_{z}^{\ast}(t)]\hat{z}(t)+[C^{\rho}(t)-f_{y}^{\ast}(t)]\hat
{y}(t)+[D^{\rho}(t)-f_{q}^{\ast}(t)]\hat{q}(t)\\
&  & +\int_{0}^{1}[f_{u}(x^{\ast}(t)+\lambda\rho(\hat{x}(t)+\tilde{x}_{\rho
}(t)),z^{\ast}(t)+\lambda\rho(\hat{z}(t)+\tilde{z}_{\rho}(t)),y^{\ast
}(t)+\lambda\rho(\hat{y}(t)+\tilde{y}_{\rho}(t)),\\
&  & q^{\ast}(t)+\lambda\rho(\hat{q}(t)+\tilde{q}_{\rho}(t)),u^{\ast
}(t)+\lambda\rho(u(t)-u^{\ast}(t)),t)-f_{u}^{\ast}(t)](u(t)-u^{\ast
}(t))d\lambda
\end{array}
\]
and
\[
\begin{array}
[c]{rrl}%
A_{1}^{\rho}(t) & = & \int_{0}^{1}g_{x}(x^{\ast}(t)+\lambda\rho(\hat
{x}(t)+\tilde{x}_{\rho}(t)),z^{\ast}(t)+\lambda\rho(\hat{z}(t)+\tilde{z}%
_{\rho}(t)),y^{\ast}(t)+\lambda\rho(\hat{y}(t)+\tilde{y}_{\rho}(t)),\\
&  & q^{\ast}(t)+\lambda\rho(\hat{q}(t)+\tilde{q}_{\rho}(t)),u^{\ast
}(t)+\lambda\rho(u(t)-u^{\ast}(t)),t)d\lambda,\\
B_{1}^{\rho}(t) & = & \int_{0}^{1}g_{z}(x^{\ast}(t)+\lambda\rho(\hat
{x}(t)+\tilde{x}_{\rho}(t)),z^{\ast}(t)+\lambda\rho(\hat{z}(t)+\tilde{z}%
_{\rho}(t)),y^{\ast}(t)+\lambda\rho(\hat{y}(t)+\tilde{y}_{\rho}(t)),\\
&  & q^{\ast}(t)+\lambda\rho(\hat{q}(t)+\tilde{q}_{\rho}(t)),u^{\ast
}(t)+\lambda\rho(u(t)-u^{\ast}(t)),t)d\lambda,\\
C_{1}^{\rho}(t) & = & \int_{0}^{1}g_{y}(x^{\ast}(t)+\lambda\rho(\hat
{x}(t)+\tilde{x}_{\rho}(t)),z^{\ast}(t)+\lambda\rho(\hat{z}(t)+\tilde{z}%
_{\rho}(t)),y^{\ast}(t)+\lambda\rho(\hat{y}(t)+\tilde{y}_{\rho}(t)),\\
&  & q^{\ast}(t)+\lambda\rho(\hat{q}(t)+\tilde{q}_{\rho}(t)),u^{\ast
}(t)+\lambda\rho(u(t)-u^{\ast}(t)),t)d\lambda,
\end{array}
\]
\[%
\begin{array}
[c]{rrl}%
D_{1}^{\rho}(t) & = & \int_{0}^{1}g_{q}(x^{\ast}(t)+\lambda\rho(\hat
{x}(t)+\tilde{x}_{\rho}(t)),z^{\ast}(t)+\lambda\rho(\hat{z}(t)+\tilde{z}%
_{\rho}(t)),y^{\ast}(t)+\lambda\rho(\hat{y}(t)+\tilde{y}_{\rho}(t)),\\
&  & q^{\ast}(t)+\lambda\rho(\hat{q}(t)+\tilde{q}_{\rho}(t)),u^{\ast
}(t)+\lambda\rho(u(t)-u^{\ast}(t)),t)d\lambda,\\
E_{1}^{\rho}(t) & = & [A_{1}^{\rho}(t)-g_{x}^{\ast}(t)]\hat{x}(t)+[B_{1}%
^{\rho}(t)-g_{z}^{\ast}(t)]\hat{z}(t)+[C_{1}^{\rho}(t)-g_{y}^{\ast}(t)]\hat
{y}(t)+[D_{1}^{\rho}(t)-g_{q}^{\ast}(t)]\hat{q}(t)\\
&  & +\int_{0}^{1}[g_{u}(x^{\ast}(t)+\lambda\rho(\hat{x}(t)+\tilde{x}_{\rho
}(t)),z^{\ast}(t)+\lambda\rho(\hat{z}(t)+\tilde{z}_{\rho}(t)),y^{\ast
}(t)+\lambda\rho(\hat{y}(t)+\tilde{y}_{\rho}(t)),\\
&  & q^{\ast}(t)+\lambda\rho(\hat{q}(t)+\tilde{q}_{\rho}(t)),u^{\ast
}(t)+\lambda\rho(u(t)-u^{\ast}(t)),t)-g_{u}^{\ast}(t)](u(t)-u^{\ast
}(t))d\lambda.
\end{array}
\]
Thus%
\[
\left\{
\begin{array}
[c]{rrl}%
-d\tilde{y}_{\rho}(t) & = & \left(  A^{\rho}(t)\tilde{x}_{\rho}(t)+B^{\rho
}(t)\tilde{z}_{\rho}(t)+C^{\rho}(t)\tilde{y}_{\rho}(t)+D^{\rho}(t)\tilde
{q}_{\rho}(t)+E^{\rho}(t)\right)  dt\\
&  & +\left(  A_{1}^{\rho}(t)\tilde{x}_{\rho}(t)+B_{1}^{\rho}(t)\tilde
{z}_{\rho}(t)+C_{1}^{\rho}(t)\tilde{y}_{\rho}(t)+D_{1}^{\rho}(t)\tilde
{q}_{\rho}(t)+E_{1}^{\rho}(t)\right)  dB_{t} -\tilde{q}_{\rho}(t)dW(t),\\
\tilde{y}_{\rho}(T) & = & 0.
\end{array}
\right.
\]
Using Lemma 2.1 to $|\tilde{y}_{\rho}(t)|^{2}$, we get%
\[%
\begin{array}
[c]{rl}
& E|\tilde{y}_{\rho}(t)|^{2}+E\int_{t}^{T}\|\tilde{q}_{\rho}(s)\|^{2}ds\\
= & 2E\int_{t}^{T}(\tilde{y}_{\rho}(t),A^{\rho}(s)\tilde{x}_{\rho}(s)+B^{\rho
}(s)\tilde{z}_{\rho}(s)+C^{\rho}(s)\tilde{y}_{\rho}(s)+D^{\rho}(s)\tilde
{q}_{\rho}(s)+E^{\rho}(s))ds\\
& +E\int_{t}^{T}\|A_{1}^{\rho}(s)\tilde{x}_{\rho}(s)+B_{1}^{\rho}(s)\tilde
{z}_{\rho}(s)+C_{1}^{\rho}(s)\tilde{y}_{\rho}(s)+D_{1}^{\rho}(s)\tilde
{q}_{\rho}(s)+E_{1}^{\rho}(s)\|^{2}ds\\
\leq & K(E\int_{t}^{T}|\tilde{y}_{\rho}(s)|^{2}ds+E\int_{t}^{T}|\tilde
{x}_{\rho}(s)|^{2}ds)+\beta_{1}(E\int_{t}^{T}\|\tilde{z}_{\rho}(s)\|^{2}%
ds+E\int_{t}^{T}\|\tilde{q}_{\rho}(s)\|^{2}ds+J_{\rho}),
\end{array}
\]
where $K>0,\ 0<\beta_{1}<\frac{1}{2}$ are constants and 
$
J_{\rho}  =  E\int_{t}^{T}|E^{\rho}(s)|^{2}ds+E\int_{t}^{T}|E_{1}^{\rho
}(s)|^{2}ds.
$

Similar analysis shows that
\[
\left\{
\begin{array}
[c]{llll}
-d\tilde{x}_{\rho}(t)& =& \left(  A^{\prime\rho}(t)\tilde{x}_{\rho}
(t)+B^{\prime\rho}(t)\tilde{z}_{\rho}(t)+C^{\prime\rho}(t)\tilde{y}_{\rho
}(t)+D^{\prime\rho}(t)\tilde{q}_{\rho}(t)+E^{\prime\rho}(t)\right)  dt 
\\
&&+\left(  A_{1}^{\prime\rho}(t)\tilde{x}_{\rho}(t)+B_{1}
^{\prime\rho}(t)\tilde{z}_{\rho}(t)+C_{1}^{\prime\rho}(t)\tilde{y}_{\rho
}(t)+D_{1} ^{\prime\rho}(t)\tilde{q}_{\rho}(t)+E_{1}^{\prime\rho}(t)\right)
dW_{t}-\tilde{z}_{\rho}(t)dB(t),  \\
\tilde{x}_{\rho}(0)& = &0, 
\end{array}
\right.
\]
where $H^{\prime\rho}(t)$ and $H_{1}^{\prime\rho}(t)$ ($H=A,$ $B,$ $C$, $D$
and $E$) are similarly defined as above.

It yields that%
\[%
\begin{array}
[c]{rl}
& E|\tilde{x}_{\rho}(t)|^{2}+E\int_{t}^{T}\|\tilde{z}_{\rho}(s)\|^{2}ds\\
= & 2E\int_{t}^{T}(\tilde{x}_{\rho}(t),A^{\prime\rho}(s)\tilde{x}_{\rho
}(s)+B^{\prime\rho}(s)\tilde{z}_{\rho}(s)+C^{\prime\rho}(s)\tilde{y}_{\rho
}(s)+D^{\prime\rho}(s)\tilde{q}_{\rho}(ts)+E^{\prime\rho}(s))ds\\
& +E\int_{t}^{T}\|A_{1}^{\prime\rho}(s)\tilde{x}_{\rho}(s)+B_{1}^{\prime\rho
}(s)\tilde{z}_{\rho}(s)+C_{1}^{\prime\rho}(s)\tilde{y}_{\rho}(s)+D_{1}%
^{\prime\rho}(s)\tilde{q}_{\rho}(s)+E_{1}^{\prime\rho}(s)\|^{2}ds\\
\leq & K(E[\int_{t}^{T}|\tilde{y}_{\rho}(s)|^{2}ds+\int_{t}^{T}|\tilde
{x}_{\rho}(s)|^{2}ds]+\beta_{2}E[\int_{t}^{T}\|\tilde{z}_{\rho}(s)\|^{2}%
ds+\int_{t}^{T}\|\tilde{q}_{\rho}(s)\|^{2}ds+J_{\rho}^{\prime}]
\end{array}
\]
where $K>0,$ $0<\beta_{2}<\frac{1}{2}$ are constants and%
$
J_{\rho}^{\prime} =  E[\int_{t}^{T}|E^{\prime\rho}(s)|^{2}ds+\int_{t}%
^{T}|E_{1}^{\prime\rho}(s)|^{2}ds].
$
Since $0<\beta_{1}+\beta_{2}<1$ and $K>0$, there exists $K_{1}>0$ such that%
\[%
\begin{array}
[c]{rl}
& E|\tilde{x}_{\rho}(t)|^{2}+E|\tilde{y}_{\rho}(t)|^{2}+(1-\beta_{1}-\beta
_{2})(E\int_{t}^{T}\|\tilde{z}_{\rho}(s)\|^{2}ds+E\int_{t}^{T}\|\tilde
{q}_{\rho}(s)\|^{2}ds)\\
\leq & K_{1}(E\int_{t}^{T}|\tilde{x}_{\rho}(s)|^{2}ds+E\int_{t}^{T}|\tilde
{y}_{\rho}(s)|^{2}ds+J_{\rho}+J_{\rho}^{\prime}).
\end{array}
\]
Since the Lebesgue dominated convergence theorem implies
$
\lim\limits_{\rho\rightarrow0}J_{\rho}=0,\ \lim\limits_{\rho\rightarrow
0}J_{\rho}^{\prime}=0,
$
we obtain the result by Gronwall's inequality.
\end{proof}

\subsubsection{Variational inequality}

In this subsection, we apply Ekeland's variational principle \cite{Ekeland} to
deal with initial-terminal state constraints
\[
E(\psi(x_{T}^{(\xi,\eta,u(\cdot))}))=a,\ E(h(y_{0}^{(\xi,\eta,u(\cdot))}))=b.
\]
Define
\begin{equation}%
\begin{array}
[c]{rrl}%
F_{\varepsilon}((\xi,\eta,u(\cdot))) & \triangleq & \{|E(\psi(x_{T}^{(\xi
,\eta,u(\cdot))}))-a|^{2}+|E(h(y_{0}^{(\xi,\eta,u(\cdot))}))-b|^{2}\\
&  & +(\max(0,E[\phi(x^{\ast}(T))-\phi(x_{T}^{(\xi,\eta,u(\cdot))}%
)]+\varepsilon))^{2}\\
&  & +(\max(0,E[\gamma(y^{\ast}(0))-\gamma(y_{0}^{(\xi,\eta,u(\cdot
))})]+\varepsilon))^{2}\}^{\frac{1}{2}},
\end{array}
\tag{2.9}%
\end{equation}
where $a$ and $b$ are the given initial and terminal state constraints and
$\varepsilon$ is an arbitrary positive constant.

It is easy to check that the mapping $|E(\psi(x_{T}^{(\xi,\eta,u(\cdot
))}))-a|^{2},$ $|E(h(y_{0}^{(\xi,\eta,u(\cdot))}))-b|^{2},$ $\phi(x_{T}%
^{(\xi,\eta,u(\cdot))})$ and $\gamma(y_{0}^{(\xi,\eta,u(\cdot))})$ are all
continuous functionals from $U$ to $%
\mathbb{R}
$.

\begin{theorem}
Suppose (H1)$\sim$(H3). Let $(\xi^{\ast},\eta^{\ast},u^{\ast}(\cdot))$ be an
optimal solution to (2.4). Then there exist $h_{0},h_{1},h_{2},h_{3}\in%
\mathbb{R}
$ with $h_{0},h_{1}\leq0$ and $|h_{0}|+|h_{1}|+|h_{2}|+|h_{3}|\neq0$ such that
the following variational inequality holds%
\begin{equation}%
\begin{array}
[c]{l}%
h_{3}E\langle\psi_{x}(x^{\ast}(T)),\hat{x}(T)\rangle+h_{2}E\langle
h_{y}(y^{\ast}(0)),\hat{y}(0)\rangle
+h_{1}E\langle\phi_{x}(x^{\ast}(T)),\hat{x}(T)\rangle+h_{0}E\langle\gamma
_{y}(y^{\ast}(0)),\hat{y}(0)\rangle\geq0
\end{array}
\tag{2.10}%
\end{equation}
\ where $\hat{x}(T)$ is the solution $\hat{x}(\cdot)$ of (2.5) at time $T$,
and $\hat{y}(0)$ is the solution $\hat{y}(\cdot)$ of (2.5) at time $0$.
\end{theorem}

\begin{proof}
It is easy to check that $F_{\varepsilon}(\cdot)$ is continuous on $U$ such
that%
\[%
\begin{array}
[c]{l}%
F_{\varepsilon}(\xi^{\ast},\eta^{\ast},u^{\ast}(\cdot))=\sqrt{2}\varepsilon;\\
F_{\varepsilon}(\xi,\eta,u(\cdot))>0,\;\forall(\xi,\eta,u(\cdot))\in U;\\
F_{\varepsilon}(\xi^{\ast},\eta^{\ast},u^{\ast}(\cdot))\leq\underset{(\xi
,\eta,u(\cdot))\in U}{\inf}F_{\varepsilon}(\xi,\eta,u(\cdot))+\sqrt
{2}\varepsilon.
\end{array}
\]
Thus, from Ekeland's variational principle \cite{Ekeland}, $\exists
(\xi^{\varepsilon},\eta^{\varepsilon},u^{\epsilon}(\cdot))\in U$ such that%
\[%
\begin{array}
[c]{l}%
\text{(i) }F_{\varepsilon}(\xi^{\varepsilon},\eta^{\varepsilon},u^{\varepsilon
}(\cdot))\leq F_{\varepsilon}(\xi^{\ast},\eta^{\ast},u^{\ast}(\cdot));\\
\text{(ii) }d((\xi^{\ast},\eta^{\ast},u^{\ast}(\cdot)),(\xi^{\varepsilon}%
,\eta^{\varepsilon},u^{\varepsilon}(\cdot)))\leq\sqrt{2\varepsilon};\\
\text{(iii) }F_{\varepsilon}(\xi,\eta,u(\cdot))+\sqrt{2\varepsilon}d((\xi
,\eta,u(\cdot)),(\xi^{\varepsilon},\eta^{\varepsilon},u^{\varepsilon}%
(\cdot)))\geq F_{\varepsilon}(\xi^{\varepsilon},\eta^{\varepsilon
},u^{\varepsilon}(\cdot)),\ \forall(\xi,\eta,u(\cdot))\in U.
\end{array}
\]
For any $(\xi,\eta,u(\cdot))\in U$ and $0\leq\rho\leq1,$ set
$
(\xi_{\rho}^{\varepsilon},\eta_{\rho}^{\varepsilon},u_{\rho}^{\varepsilon
}(\cdot)) =(\xi^{\varepsilon}+\rho(\xi-\xi^{\varepsilon}),\eta^{\varepsilon
}+\rho(\eta-\eta^{\varepsilon}),u^{\varepsilon}(\cdot)+\rho(u(\cdot
)-u^{\varepsilon}(\cdot))).
$
Let $(x_{\rho}^{\varepsilon}(\cdot),y_{\rho}^{\varepsilon}(\cdot),z_{\rho
}^{\varepsilon}(\cdot),q_{\rho}^{\varepsilon}(\cdot))$ (resp. $(x^{\varepsilon
}(\cdot),y^{\varepsilon}(\cdot),z^{\varepsilon}(\cdot),q^{\varepsilon}%
(\cdot))$) be the solution of (2.2) under $(\xi_{\rho}^{\varepsilon}%
,\eta_{\rho}^{\varepsilon},u_{\rho}^{\varepsilon}(\cdot))$ (resp.$(\xi
^{\varepsilon},\eta^{\varepsilon},u^{\varepsilon}(\cdot)$)$,$ and $(\hat
{x}^{\varepsilon}(\cdot),\hat{y}^{\varepsilon}(\cdot),\hat{z}^{\varepsilon
}(\cdot),\hat{q}^{\varepsilon}(\cdot))$ be the solution of (2.5) in which
$(\xi^{\ast},\eta^{\ast},u^{\ast}(\cdot))$ is substituted by $(\xi
^{\varepsilon},\eta^{\varepsilon},u^{\varepsilon}(\cdot))$.

From (iii), we know that%
\begin{equation}
F_{\varepsilon}(\xi_{\rho}^{\varepsilon},\eta_{\rho}^{\varepsilon},u_{\rho
}^{\varepsilon}(\cdot))-F_{\varepsilon}(\xi^{\varepsilon},\eta^{\varepsilon
},u^{\varepsilon}(\cdot))+\sqrt{2\varepsilon}d((\xi_{\rho}^{\varepsilon}%
,\eta_{\rho}^{\varepsilon},u_{\rho}^{\varepsilon}(\cdot)),(\xi^{\varepsilon
},\eta^{\varepsilon},u^{\varepsilon}(\cdot)))\geq0. \tag{2.11}%
\end{equation}

On the other hand, similarly to Lemma 2.4 we have%
\begin{align*}
\lim\limits_{\rho\rightarrow0}\sup\limits_{0\leq t\leq T}E[\rho^{-1}[x_{\rho
}^{\varepsilon}(t)-x^{\varepsilon}(t)]-\hat{x}^{\varepsilon}(t)]^{2}  &  =0,\\
\lim\limits_{\rho\rightarrow0}\sup\limits_{0\leq t\leq T}E[\rho^{-1}[y_{\rho
}^{\varepsilon}(t)-y^{\varepsilon}(t)]-\hat{y}^{\varepsilon}(t)]^{2}  &  =0.
\end{align*}

This leads to the following expansions%
\[%
\begin{array}
[c]{lll}
&  & |E(\phi(x_{\rho}^{\varepsilon}(T))-a|^{2}-|E(\phi(x^{\varepsilon
}(T))-a|^{2}
=  2\rho\langle E(\phi(x^{\varepsilon}(T))-a,E[\phi_{x}(x^{\varepsilon
}(T))\hat{x}^{\varepsilon}(T)]\rangle+o(\rho),\\
&  & |E(h(y_{\rho}^{\varepsilon}(0))-b|^{2}-|E(h(y{\varepsilon}(0))-b|^{2}= 2\rho\langle E(h(y^{\varepsilon}(0))-b,E[h_{y}(y^{\varepsilon}%
(0))\hat{y}^{\varepsilon}(0)]\rangle+o(\rho).
\end{array}
\]

Applying the linearization technique, then%
\[%
\begin{array}
[c]{lll}%
E[\phi(x_{\rho}^{\varepsilon}(T))-\phi(x^{\varepsilon}(T))] = E[\int_{0}%
^{T}\phi_{x}(x^{(\xi^{\varepsilon}+\lambda\rho\hat{\xi}^{\varepsilon}%
,\eta^{\varepsilon} +\lambda\rho\hat{\eta}^{\varepsilon},u^{\varepsilon}%
(\cdot)+\lambda\rho\hat{u}^{\varepsilon})}(T))d\lambda\cdot\rho\hat
{x}^{\varepsilon}(T)]; &  & \\
E[\gamma(y_{\rho}^{\varepsilon}(0))-\gamma(y^{\varepsilon}(0))] = E[\int%
_{0}^{T}\gamma_{y}(y ^{(\xi^{\varepsilon}+\lambda\rho\hat{\xi}^{\varepsilon
},\eta^{\varepsilon}+\lambda\rho\hat{\eta}^{\varepsilon},u^{\varepsilon}%
(\cdot)+\lambda\rho\hat{u}^{\varepsilon}(\cdot))}(0))d\lambda\cdot\rho\hat
{y}^{\varepsilon}(0)]. &  &
\end{array}
\]
So we have the following expansions%
\[%
\begin{array}
[c]{rl}
& (E[\phi(x^{\ast}(T))-\phi(x_{\rho}^{\varepsilon}(T)]+\varepsilon
)^{2}-(E[\phi(x^{\ast}(T))-\phi(x^{\varepsilon}(T)]+\varepsilon)^{2}\\
= & -2\rho E[\phi_{x}(x^{\varepsilon}(T))\cdot\rho\hat{x}_{T}^{\varepsilon
}][E[\phi(x^{\ast}(T))-\phi(x^{\varepsilon}(T)]+\varepsilon],\\
& (E[\gamma(y^{\ast}(0)-\gamma(y_{\rho}^{\varepsilon}(0)]+\varepsilon
)^{2}-(E[\gamma(y^{\ast} (0))-\gamma(y^{\varepsilon}(0))]+\varepsilon)^{2}\\
= & -2\rho E[\gamma_{y}(y^{\varepsilon}(0))\cdot\rho\hat{y}_{0}^{\varepsilon
}][E[\gamma(y^{\ast}(0))-\gamma(y^{\varepsilon}(0))]+\varepsilon].
\end{array}
\]
For the given $\varepsilon$, we consider the following four cases: \newline
Case 1.  There exists $\rho_{0}>0$ such that%
\[
E[\phi(x_{\rho}^{\varepsilon}(T))-\phi(x^{\varepsilon}(T))]+\varepsilon
>0,\ E[\gamma(y_{\rho}^{\varepsilon}(0))-\gamma(y^{\varepsilon}%
(0))]+\varepsilon>0
\]
for all $\rho\in(0,\rho_{0})$.

In this case,%
\[%
\begin{array}
[c]{rl}
& \lim\limits_{\rho\rightarrow0}\frac{F_{\varepsilon}(\xi_{\rho}^{\varepsilon
},\eta_{\rho}^{\varepsilon},u_{\rho}^{\varepsilon}(\cdot))-F_{\varepsilon}%
(\xi^{\varepsilon},\eta^{\varepsilon},u^{\varepsilon}(\cdot))}{\rho}\\
= & \lim\limits_{\rho\rightarrow0}\frac{1}{F_{\varepsilon}(\xi_{\rho
}^{\varepsilon},\eta_{\rho}^{\varepsilon},u_{\rho}^{\varepsilon}%
(\cdot))+F_{\varepsilon}(\xi^{\varepsilon},\eta^{\varepsilon},u^{\varepsilon
}(\cdot))}\frac{F_{\varepsilon}^{2}(\xi_{\rho}^{\varepsilon},\eta_{\rho
}^{\varepsilon},u_{\rho}^{\varepsilon}(\cdot))-F_{\varepsilon}^{2}%
(\xi^{\varepsilon},\eta^{\varepsilon},u^{\varepsilon}(\cdot))}{\rho}\\
= & \frac{1}{F_{\varepsilon}(\xi^{\varepsilon},\eta^{\varepsilon
},u^{\varepsilon}(\cdot))}\{\langle E\psi(x^{\varepsilon}(T))-a,E[\psi
_{x}(x^{\varepsilon}(T))\hat{x}^{\varepsilon}(T)]\rangle+\langle Eh(y^{\varepsilon}(0))-b,E[h_{y}(y^{\varepsilon}(0))\hat
{y}^{\varepsilon}(0)]\rangle\\
& -\langle E[\phi(x^{\ast}(T))-\phi(x^{\varepsilon}(T))]+\varepsilon,
E[\phi_{x}(x^{\varepsilon}(T))\hat{x}^{\varepsilon}(T)]\rangle\\
& -\langle E[\gamma(y^{\ast}(0) )-\gamma(y^{\varepsilon}(0))]+\varepsilon,
E[\gamma_{y}(y^{\varepsilon}(0))\hat{y}^{\varepsilon}(0)]\rangle\}.
\end{array}
\]
Dividing (2.11) by $\rho$ and sending $\rho$ to $0$, we obtain%
\[%
\begin{array}
[c]{rl}
& h_{3}^{\varepsilon}E\langle\psi_{x}(x^{\varepsilon}(T)),\hat{x}%
^{\varepsilon}(T)\rangle+h_{2}^{\varepsilon}E\langle h_{y}(y^{\varepsilon
}(0)),\hat{y}^{\varepsilon}(0)\rangle\\
& +h_{1}^{\varepsilon}E\langle\phi_{x}(x^{\varepsilon}(T)),\hat{x}%
^{\varepsilon}(T)\rangle+h_{0}^{\varepsilon}E\langle\gamma_{y}(y^{\varepsilon
}(0),\hat{y}^{\varepsilon}(0)\rangle
\geq  -\sqrt{2\varepsilon}d((\xi^{\varepsilon},\eta^{\varepsilon
},u^{\varepsilon}(\cdot)),(\xi,\eta,u(\cdot))),
\end{array}
\]
where%
\begin{align*}
h_{0}^{\varepsilon}  &  =-\frac{1}{F_{\varepsilon}(\xi^{\varepsilon}%
,\eta^{\varepsilon},u^{\varepsilon}(\cdot))}[E[\gamma(y^{\ast}(0))-\gamma
(y^{\varepsilon}(0))]+\varepsilon]\leq0,\\
h_{1}^{\varepsilon}  &  =-\frac{1}{F_{\varepsilon}(\xi^{\varepsilon}%
,\eta^{\varepsilon},u^{\varepsilon}(\cdot))}[E[\phi(x^{\ast}(T))-\phi
(x^{\varepsilon}(T)]+\varepsilon]\leq0,\\
h_{2}^{\varepsilon}  &  =\frac{1}{F_{\varepsilon}(\xi^{\varepsilon}%
,\eta^{\varepsilon},u^{\varepsilon}(\cdot))}[Eh(y^{\varepsilon}(0))-b],\\
h_{3}^{\varepsilon}  &  =\frac{1}{F_{\varepsilon}(\xi^{\varepsilon}%
,\eta^{\varepsilon},u^{\varepsilon}(\cdot))}[E\psi(x^{\varepsilon}(T))-a].
\end{align*}
Case 2. There exists a position sequence \{$\rho_{n}$\} satisfying $\rho
_{n}\rightarrow0$ such that%
\[%
\begin{array}
[c]{lll}
&  & E[\phi(x_{\rho_{n}}^{\varepsilon}(T))-\phi(x^{\ast}(T))]+\varepsilon
\leq0,\ E[\gamma(y_{\rho_{n}}^{\varepsilon}(0))-\gamma(y^{\ast}(0))]+\varepsilon>
0.
\end{array}
\]
Then
\[%
\begin{array}
[c]{rl}%
F_{\varepsilon}(\xi_{\rho_{n}}^{\varepsilon},\eta_{\rho_{n}}^{\varepsilon
},u_{\rho_{n}}^{\varepsilon}(\cdot)) = & \{|E\psi(x_{\rho_{n}}^{\varepsilon
}(T))-a|^{2}+|Eh(y_{\rho_{n}}^{\varepsilon}(0))-b|^{2}\\
& +(\max(0,E[\gamma(y^{\ast}(0))-\gamma(y_{\rho_{n}}^{\varepsilon
}(0))]+\varepsilon))^{2}\}^{\frac{1}{2}}.
\end{array}
\]
For sufficiently large $n$, since $F_{\varepsilon}(\cdot)$ is continuous, we
conclude%
\[%
\begin{array}
[c]{cl}%
F_{\varepsilon}(\xi^{\varepsilon},\eta^{\varepsilon},u^{\varepsilon}(\cdot))
= & \{|E\psi(x^{\varepsilon}(T))-a|^{2}\}^{\frac{1}{2}}+|Eh(y^{\varepsilon
}(0))-b|^{2}\\
& +(\max(0,E[\gamma(y^{\ast}(0))-\gamma(y^{\varepsilon}(0))]+\varepsilon
))^{2}\}^{\frac{1}{2}}.
\end{array}
\]
Now%
\[%
\begin{array}
[c]{rl}
& \lim\limits_{n\rightarrow\infty}\frac{F_{\varepsilon}(\xi_{\rho_{n}%
}^{\varepsilon},\eta_{\rho_{n}}^{\varepsilon},u_{\rho_{n}}^{\varepsilon}%
(\cdot))-F_{\varepsilon}(\xi^{\varepsilon},\eta^{\varepsilon},u^{\varepsilon
}(\cdot))}{\rho_{n}}\\
= & \lim\limits_{n\rightarrow\infty}\frac{1}{F_{\varepsilon}(\xi_{\rho_{n}%
}^{\varepsilon},\eta_{\rho_{n}}^{\varepsilon},u_{\rho_{n}}^{\varepsilon}%
(\cdot))+F_{\varepsilon}(\xi^{\varepsilon},\eta^{\varepsilon},u^{\varepsilon
}(\cdot))}\frac{F_{\varepsilon}^{2}(\xi_{\rho_{n}}^{\varepsilon},\eta
_{\rho_{n}}^{\varepsilon},u_{\rho_{n}}^{\varepsilon}(\cdot)))-F_{\varepsilon
}^{2}(\xi^{\varepsilon},\eta^{\varepsilon},u^{\varepsilon}(\cdot))}{\rho_{n}%
}\\
= & \frac{1}{F_{\varepsilon}(\xi^{\varepsilon},\eta^{\varepsilon
},u^{\varepsilon}(\cdot))}\{\langle E(\psi(x^{\varepsilon}(T))-a,E[\psi
_{x}(x^{\varepsilon}(T))\hat{x}^{\varepsilon}(T)]\rangle\\
& +\langle Eh(y^{\varepsilon}(0))-b, E[h_{y}(y^{\varepsilon}(0))\hat
{y}^{\varepsilon}(0)]\rangle\\
& -\langle E[\gamma(y^{\ast}(0))-\gamma(y^{\varepsilon}(0))]+\varepsilon,
E[\gamma_{y}(y^{\varepsilon}(0))\hat{y}^{\varepsilon}(0)]\rangle\}.
\end{array}
\]
Similar to Case 1 we get%
\[%
\begin{array}
[c]{rl}
& h_{3}^{\varepsilon}E\langle\psi_{x}(x^{\varepsilon}(T)),\hat{x}%
^{\varepsilon}(T)\rangle+h_{2}^{\varepsilon}E\langle h_{y}(y^{\varepsilon
}(0)),\hat{y}^{\varepsilon}(0)\rangle+h_{0}^{\varepsilon}E\langle\gamma
_{y}(y^{\varepsilon}(0)),\hat{y}^{\varepsilon}(0)\rangle\\
\geq & -\sqrt{2\varepsilon}d((\xi^{\varepsilon},\eta^{\varepsilon
},u^{\varepsilon}(\cdot)),(\xi,\eta,u(\cdot))),
\end{array}
\]
where%
\begin{align*}
h_{0}^{\varepsilon}  &  =-\frac{1}{F_{\varepsilon}(\xi^{\varepsilon}%
,\eta^{\varepsilon},u^{\varepsilon}(\cdot))}[E[\gamma(y^{\ast}(0))-\gamma
(y^{\varepsilon}(0))]+\varepsilon]\leq0,\\
h_{1}^{\varepsilon}  &  =0,\\
h_{2}^{\varepsilon}  &  =\frac{1}{F_{\varepsilon}(\xi^{\varepsilon}%
,\eta^{\varepsilon},u^{\varepsilon}(\cdot))}[E(h(y^{\varepsilon}(0))-b],\\
h_{3}^{\varepsilon}  &  =\frac{1}{F_{\varepsilon}(\xi^{\varepsilon}%
,\eta^{\varepsilon},u^{\varepsilon}(\cdot))}[E(\psi(x^{\varepsilon}(T))-a].
\end{align*}
Case 3. There exists a positive sequence \{$\rho_{n}$\} satisfying $\rho
_{n}\rightarrow0$ such that%
\[%
\begin{array}
[c]{lll}%
E[\phi(x_{\rho_{n}}^{\varepsilon}(T))-\phi(x^{\ast}(T))]+\varepsilon> 0, \
E[\gamma(y_{\rho_{n}}^{\varepsilon}(0))-\gamma(y^{\ast}(0))]+\varepsilon
\leq0. &  &
\end{array}
\]
Case 4. There exists a positive sequence \{$\rho_{n}$\} satisfying $\rho
_{n}\rightarrow0$ such that%
\[%
\begin{array}
[c]{lll}%
E[\phi(x_{\rho_{n}}^{\varepsilon}(T))-\phi(x^{\ast}(T))]+\varepsilon\leq0,  \
E[\gamma(y_{\rho_{n}}^{\varepsilon}(0))-\gamma(y^{\ast}(0))]+\varepsilon
\leq0. &  &
\end{array}
\]
Similar techniques can be used to both Case 3 and Case 4.

In summary, for all those cases, we have $h_{0}^{\varepsilon}\leq
0,h_{1}^{\varepsilon}\leq0$ and $|h_{0}^{\varepsilon}|^{2}+|h_{1}%
^{\varepsilon}|^{2}+|h_{2}^{\varepsilon}|^{2}+|h_{3}^{\varepsilon}|^{2}=1$ by
the definition of $F_{\varepsilon}(\cdot)$. Then there exists a convergent
subsequence of $( h_{3}^{\varepsilon},h_{2}^{\varepsilon},h_{1}^{\varepsilon
},h_{0}^{\varepsilon})$ \ whose limit is denoted by $( h_{3},h_{2},h_{1}%
,h_{0})$. On the other hand, it is easy to check that $\hat{x}^{\varepsilon
}(T)\rightarrow0,\hat{y}^{\epsilon}(0) \rightarrow0$ as $\varepsilon
\rightarrow0.$ Thus (2.10) holds.
\end{proof}

\subsubsection{Maximum principle}

In this subsection we derive the maximum principle for the case where
$l(x,z,y,q,u,t)=0,\ \chi(x)=0,\ \lambda(y)=0$ and then present the results for
the general case in subsection 2.3.4. To this end, we introduce the adjoint
process $(m(\cdot),p(\cdot))$ and $(n(\cdot),\delta(\cdot))$ associated with
the optimal solution $(\xi^{\ast},\eta^{\ast},u^{\ast}(\cdot))$ to (2.2),
which is the solution of the following time-symmetric FBDSDE
\begin{equation}
\left\{
\begin{array}
[c]{rrl}%
-dm(t) & = & -(F_{x}^{\ast}(t)m(t)+G_{x}^{\ast}(t)p(t)+f_{x}^{\ast
}(t)n(t)+g_{x}^{\ast}(t)\delta(t))dt\\
&  & +(F_{z}^{\ast}(t)m(t)+G_{z}^{\ast}(t)p(t)+f_{z}^{\ast}(t)n(t)+g_{z}
^{\ast}(t)\delta(t))dB_{t} -p_{t}dW_{t},\\
m(T) & = & -(h_{3}\psi_{x}(x_{T}^{\ast})+h_{1}\phi_{x}(x_{T}^{\ast})),\\
dn(t) & = & (F_{y}^{\ast}(t)m(t)+G_{y}^{\ast}(t)p(t)+f_{y}^{\ast}%
(t)n(t)+g_{y}^{\ast}(t)\delta(t))dt\\
&  & +(F_{q}^{\ast}(t)m(t)+G_{q}^{\ast}(t)p(t)+f_{q}^{\ast}(t)n(t)+g_{q}%
^{\ast}(t)\delta(t))dW_{t}-\delta_{t}dB_{t},\\
n(0) & = & h_{2}h_{y}(y_{0}^{\ast})+h_{0}\gamma_{y}(y_{0}^{\ast}),
\end{array}
\right.  \tag{2.12}%
\end{equation}
where $F_{k}^{\ast}(t),$ $f_{k}^{\ast}(t),G_{k}^{\ast}(t),g_{k}^{\ast}(t)$ for
$k=x,y,z,q$ are defined as in $(2.5)$. It is easy to check that there exist
unique processes $(m(\cdot),p(\cdot))$, $(n(\cdot),\delta(\cdot))$ which solve
the above equations.

\begin{theorem}
We assume (H1)$\sim$(H4). Let $(\xi^{\ast},\eta^{\ast},u^{\ast}(\cdot))$ be
optimal and $(x^{\ast}(\cdot),$ $z^{\ast}(\cdot),y^{\ast}(\cdot),q^{\ast
}(\cdot))$ be the corresponding optimal trajectory. Then for arbitrary
$(\xi,\eta,u(\cdot))$ $\in U$ we have for any $t\in[0,T]$,
\begin{equation}%
\begin{array}
[c]{l}%
\langle m(0),\xi-\xi^{\ast})\rangle\leq0,\\
\langle n(T),\eta-\eta^{\ast}\rangle\geq0,\\
\langle m(t)F_{u}^{\ast}(t)+p(t)G_{u}^{\ast}(t)+n(t)f_{u}^{\ast}%
(t)+g_{u}^{\ast}(t)\delta(t),u(t)-u^{\ast}(t) \rangle\geq0.
\end{array}
\tag{2.13}%
\end{equation}

\end{theorem}

\begin{proof}
For any $(\xi,\eta,u(\cdot))\in U$, let $(\hat{x}(\cdot),\hat{z}(\cdot
),\hat{y}(\cdot),\hat{q}(\cdot))$ be the solution to (2.5). Applying Lemma 2.1
to $\langle m(t),\hat{x}(t)\rangle+\langle n(t),\hat{y}(t)\rangle$ , we have%
\[%
\begin{array}
[c]{rl}
& E\langle m(T),\hat{x}(T)\rangle+E\langle n(T),\hat{y}(T)\rangle-E\langle
m(0),\hat{x}(0)\rangle-E\langle n(0),\hat{y}(0)\rangle\\
= & E\int_{0}^{T}m_{t}(-F_{x}^{\ast}(t)\hat{x}(t)-F_{z}^{\ast}(t)\hat
{z}(t)-F_{y}^{\ast}(t)\hat{y}(t)-F_{q}^{\ast}(t)\hat{q}(t)-F_{u}^{\ast}%
(t)\hat{u}(t))dt\\
& +E\int_{0}^{T}\hat{x}(t)(F_{x}^{\ast}(t)m(t)+G_{x}^{\ast}(t)p(t)+f_{x}%
^{\ast}(t)n(t)+g_{x}^{\ast}(t)\delta(t))dt\\
& +E\int_{0}^{T}\hat{z}(t)(F_{z}^{\ast}(t)m(t)+G_{z}^{\ast}(t)p(t)+f_{z}%
^{\ast}(t)n(t)+g_{z}^{\ast}(t)\delta(t))dt\\
& -E\int_{0}^{T}p(t)(G_{x}^{\ast}(t)\hat{x}(t)+G_{z}^{\ast}(t)\hat{z}%
(t)+G_{y}^{\ast}(t)\hat{y}(t)+G_{q}^{\ast}(t)\hat{q}(t)+G_{u}^{\ast}(t)\hat
{u}(t))dt\\
& -E\int_{0}^{T}n(t)(f_{x}^{\ast}(t)\hat{x}(t)+f_{z}^{\ast}(t)\hat{z}%
(t)+f_{y}^{\ast}(t)\hat{y}(t)+f_{q}^{\ast}(t)\hat{q}(t)+f_{u}^{\ast}(t)\hat
{u}(t))dt\\
& +E\int_{0}^{T}\hat{y}(t)(F_{y}^{\ast}(t)m(t)+G_{y}^{\ast}(t)p(t)+f_{y}
^{\ast}(t)n(t)+g_{y}^{\ast}(t)\delta(t))dt\\
& +E\int_{0}^{T}\hat{q}(t)(F_{q}^{\ast}(t)m(t)+G_{q}^{\ast}(t)p(t)+f_{q}
^{\ast}(t)n(t)+q_{q}^{\ast}(t)\delta(t))dt\\
& -E\int_{0}^{T}\delta(t)(g_{x}^{\ast}(t)\hat{x}(t)+g_{z}^{\ast}(t)\hat
{z}(t)+g_{y}^{\ast}(t)\hat{y}(t)+g_{q}^{\ast}(t)\hat{q}(t)+g_{u}^{\ast}%
(t)\hat{u}(t))dt\\
= & -E\int_{0}^{T}\langle m(t)F_{u}^{\ast}(t)+p(t)G_{u}^{\ast}(t)+n(t)f_{u}%
^{\ast}(t)+g_{u}^{\ast}(t)\delta(t),\hat{u}(t)\rangle dt.
\end{array}
\]
This yields%
\[%
\begin{array}
[c]{rl}
& E\langle n(T),\hat{y}(T)\rangle-E\langle m(0),\hat{x}(0)\rangle\\
& +E\int_{0}^{T}\langle m(t)F_{u}^{\ast}(t)+p(t)G_{u}^{\ast}(t)+n(t)f_{u}%
^{\ast}(t)+\delta(t)g_{u}^{\ast}(t),u(t)-u^{\ast}(t)\rangle dt\\
= & -E\langle m(T),\hat{x}(T)\rangle+E\langle n(0),\hat{y}(0)\rangle\\
= & h_{3}E\langle\psi_{x}(x^{\ast}(T)),\hat{x}(T)\rangle+h_{2}E\langle
h_{y}(y^{\ast}(0)),\hat{y}(0)\rangle\\
& +h_{1}E\left\langle \phi_{x}(x^{\ast}(T)),\hat{x}(T)\right\rangle
+h_{0}E\left\langle \gamma_{y}(y^{\ast}(0)),\hat{y}(0)\right\rangle \geq0.
\end{array}
\]
For any $(\xi,\eta,u(\cdot))\in U$, we get%
\[%
\begin{array}
[c]{rl}
& E\langle n(T),\eta-\eta^{\ast}\rangle-E\langle m(0),\xi-\xi^{\ast})\rangle\\
& +E\int_{0}^{T}\langle m(t)F_{u}^{\ast}(t)+p(t)G_{u}^{\ast}(t)+n(t)f_{u}%
^{\ast}(t)+\delta(t)g_{u}^{\ast}(t),u(t)-u^{\ast}(t)\rangle dt\geq0.
\end{array}
\]
Thus, it is easy to see that (2.13) holds.
\end{proof}

\subsubsection{The general case}

Define the Hamiltonian
\[
H:\mathbb{R}^{n}\times\mathbb{R}^{n\times l}\times\mathbb{R} ^{k}%
\times\mathbb{R}^{k\times d}\times\mathbb{R}^{n\times d}\times\mathbb{R}%
^{n}\times\mathbb{R}^{n\times l}\times\mathbb{R}^{k}\times\mathbb{R}^{k\times
d}\times\lbrack0,T]\rightarrow\mathbb{R}
\]
by
\[%
\begin{array}
[c]{rl}
& H(x,z,y,q,u,m,p,n,\delta,t)
=  \langle F(t,x,z,y,q,u),m(t)\rangle+\langle G(t,x,z,y,q,u),p(t)\rangle\\
& +\langle f(t,x,z,y,q,u),n(t)\rangle+\langle g(t,x,z,y,q,u),\delta
(t)\rangle+l(x,z,y,q,u,t).
\end{array}
\]
Now we consider the general case where $l(x,z,y,q,u,t)\neq0,\ \chi
(x)\neq0,\ \lambda(y)\neq0.$

Since the proof of the maximum principle is essentially similar as in the
preceding subsection, we only present the result without proof.

Let $(\xi^{\ast},\eta^{\ast},u^{\ast}(\cdot))$ be optimal to (2.4) with
$(x^{\ast}(\cdot),z^{\ast}(\cdot),y^{\ast}(\cdot),q^{\ast}(\cdot))$ being the
corresponding optimal trajectory of (2.2). We define the following adjoint
equations
\[
\left\{
\begin{array}
[c]{rrl}%
-dm(t) & = & -(F_{x}^{\ast}(t)m(t)+G_{x}^{\ast}(t)p(t)+f_{x}^{\ast
}(t)n(t)+g_{x}^{\ast}(t)\delta(t)+l_{x}^{\ast}(t))dt\\
&  & +(F_{z}^{\ast}(t)m(t)+G_{z}^{\ast}(t)p(t)+f_{z}^{\ast}(t)n(t)+g_{z}%
^{\ast}(t)\delta(t)+l_{z}^{\ast}(t))dB_{t}-p_{t}dW_{t},\\
m(T) & = & -(h_{3}\psi_{x}(x_{T}^{\ast})+h_{1}\phi_{x}(x_{T}^{\ast})),\\
dn(t) & = & (F_{y}^{\ast}(t)m(t)+G_{y}^{\ast}(t)p(t)+f_{y}^{\ast}%
(t)n(t)+g_{y}^{\ast}(t)\delta(t)+l_{y}^{\ast}(t))dt\\
&  & +(F_{q}^{\ast}(t)m(t)+G_{q}^{\ast}(t)p(t)+f_{q}^{\ast}(t)n(t)+g_{q}%
^{\ast}(t)\delta(t)+l_{q}^{\ast}(t))dW_{t}-\delta_{t}dB_{t},\\
n(0) & = & h_{2}h_{y}(y_{0}^{\ast})+h_{0}\gamma_{y}(y_{0}^{\ast}),
\end{array}
\right.
\]
where $l_{a}^{\ast}(t)=l_{a}(x^{\ast}(t),z^{\ast}(t),y^{\ast}(t),q^{\ast
}(t),u^{\ast}(t),t)$, $a=x,z,y,q$, respectively.

\begin{theorem}
We assume (H1)$\sim$(H3). Let $(\xi^{\ast},\eta^{\ast},u^{\ast}(\cdot))$ be
optimal and $(x^{\ast}(\cdot),$ $z^{\ast}(\cdot),y^{\ast}(\cdot),q^{\ast
}(\cdot))$ be the corresponding optimal trajectory. Then for arbitrary
$(\xi,\eta,u(\cdot))\in U$, we have the following inequalities hold, for any
$t\in[0,T]$,
\[
\left\{
\begin{array}
[c]{l}%
\langle m(0)-\chi_{x}(\xi^{\ast}),\xi-\xi^{\ast}\rangle\leq0,\\
\langle n(T)+\lambda_{y}(\eta^{\ast}),\eta-\eta^{\ast}\rangle\geq0,\\
\langle H_{u}(t,x^{\ast}(t),z^{\ast}(t),y^{\ast}(t),q^{\ast}(t),u^{\ast
}(t),m(t),p(t),n(t),\delta(t) ),u(t)-u^{\ast}(t)\rangle\geq0.
\end{array}
\right.
\]

\end{theorem}

\begin{remark}
Let us denoted the boundary of $K_{1}$ by $\partial K_{1}$. Set
\[
M\triangleq\{w\in\Omega|\xi^{\ast}(\omega)\in\partial K_{1}\}.
\]
Then
\[
\left\{
\begin{array}
[c]{rrl}%
m(0)-\chi_{x}(\xi^{\ast}) & \leq & 0\text{ }a.s.\text{ }on\text{ }M,\\
m(0)-\chi_{x}(\xi^{\ast}) & = & 0\text{ }a.s.\text{ }on\text{ }M^{c}%
\end{array}
\right.
\]
Similar analysis can be used to the boundaries of $K_{2}$ and $K$.
\end{remark}

\section{Applications}

In this section, we give three specific cases to illustrate the applications
of our obtained results.

\subsection{System composed of a Forward SDE and a BDSDE}

\paragraph{\bigskip\textbf{Classical formulation}}

For given $\xi\in L^{2}(\Omega,\mathcal{F}_{0},P;\mathbb{R}^{n})$ and
$u(\cdot)\in U[0,T]$, we consider the following controlled system composed of
a FSDE and a BDSDE.
\begin{equation}
\left\{
\begin{array}
[c]{rrl}%
dy(t) & = & \bar{b}(t,y(t),u(t))dt+\sigma(t,y(t),u(t))dW_t,\\
y(0) & = & b,\\
-dx(t) & = & \bar{f}(t,x(t),z(t),y(t),u(t))dt+\bar{g}(t,x(t)%
,z(t),y(t),u(t))dW_{t}-z(t)dB_t,\\
x(0) & = & \xi,
\end{array}
\right.  \tag{3.1}%
\end{equation}
where $b\in R^{k}$ is given, $x(0)=\xi\in K_{1},$ a.s, where $K_{1}$ is a
given nonempty convex subset in $R^{n}$.

Set the mappings%
\[%
\begin{array}
[c]{rl}%
\bar{b}: & \Omega\times\lbrack0,T]\times\mathbb{R}^{k}\times\mathbb{R}%
^{n\times d}\rightarrow\mathbb{R}^{k},\\
\sigma: & \Omega\times\lbrack0,T]\times\mathbb{R}^{k}\times\mathbb{R}^{n\times
d}\rightarrow\mathbb{R}^{k\times d},\\
\ \bar{f}: & \Omega\times\lbrack0,T]\times\mathbb{R} ^{n}\times\mathbb{R}%
^{n\times l}\times\mathbb{R}^{k}\times\mathbb{R}^{n\times d}\rightarrow
\mathbb{R}^{n},\\
\bar{g}: & \Omega\times\lbrack0,T]\times\mathbb{R}^{n}\times\mathbb{R}%
^{n\times l}\times\mathbb{R} ^{k}\times\mathbb{R}^{n\times d}\rightarrow
\mathbb{R}^{n\times d}.
\end{array}
\]
In this case, we regard $u(\cdot)$ and $\xi$ as the control variables. Define
the following cost function:%
\[
J(\xi,u(\cdot))=E[\int_{0}^{T}\bar{l}(t,x(t),z(t),y(t),u(t))dt+\chi
(\xi)+\lambda(y(T))+\phi(x(T))],
\]
where%
\[%
\begin{array}
[c]{lll}%
\bar{l}: \Omega\times\lbrack0,T]\times\mathbb{R}^{n}\times\mathbb{R}^{n\times
l}\times\mathbb{R}^{k}\times\mathbb{R}^{n\times d}\rightarrow\mathbb{R}, \
\chi: \mathbb{R}^{n}\rightarrow\mathbb{R},\ \lambda: \mathbb{R} ^{k}
\rightarrow\mathbb{R},\ \phi: \mathbb{R}^{n}\rightarrow\mathbb{R}. 
\end{array}
\]

We assume:

(H1) $\bar{b},\sigma,\bar{f},\bar{g},\bar{l},\chi,\lambda$ and $\phi$ are
continuous in their arguments and continuously differentiable in $(x,z,y,u)$;

(H2) There exist constants $C>0$ and $0<\alpha<\frac{1}{2}$ such that for any
$(\omega,t)\in\Omega\times\lbrack0,T],$ $u\in\mathbb{R}^{n\times d}$,
$(x_{1},z_{1},y_{1}),(x_{2},z_{2},y_{2})\in\mathbb{R}^{n}\times\mathbb{R}%
^{n\times l}\times\mathbb{R}^{k}$ the following conditions hold:
\[
|\bar{f}(t,x_{1},z_{1},y_{1},u)-\bar{f}(t,x_{2},z_{2},y_{2},u)|^{2} \leq
C(|x_{1}-x_{2}|^{2}+\|z_{1}-z_{2}\|^{2}+|y_{1}-y_{2}|^{2}),
\]
and%
\[
\|\bar{g}(t,x_{1},z_{1},y_{1},u)-\bar{g}(t,x_{2},z_{2},y_{2},u)\|^{2} \leq
C(|x_{1}-x_{2}|^{2}+|y_{1}-y_{2}|^{2})+\alpha\|z_{1}-z_{2}\|^{2}.
\]

(H3) The derivatives of $\bar{b},\sigma,\bar{f},\bar{g}$ in $(x,y,z,u)$ are
bounded; the derivatives of $\bar{l}$ in $(x,y,z,u)$ are bounded by
$C(1+|x|+|y|+|z|+\|u\|)$; the derivatives of $\chi$ and $\phi$ in $x$ are
bounded by $C(1+|x|)$; the derivatives of $\lambda$ in $y$ are bounded by
$C(1+|y|)$.

Then, for given $\xi\in L^{2}(\Omega,\mathcal{F}_{0},P;\mathbb{R}%
^{n})$ and $u(\cdot)\in U[0,T]$, there exists a unique triple
\[
(x(\cdot),y(\cdot),z(\cdot))\in M^{2}(0,T; \mathbb{R}^{n})\times M^{2}(0,T;
\mathbb{R}^{k})\times M^{2}(0,T; \mathbb{R}^{n\times l})
\]
which solves (3.1).

We assume an additional terminal state constraint $y(T)=\eta\in K_{2}$,
$a.s$., where $K_{2}$ is a given nonempty convex subset in $\mathbb{R}^{k}$.
Our stochastic control problem is%
\begin{align}
&  \inf\ \ J(\xi,u(\cdot))\tag{3.2}\\
\text{subject to }u(\cdot)  &  \in U[0,T]\text{; }\xi\in K_{1},a.s\text{,
}\eta\in K_{2}\text{, }a.s.\nonumber
\end{align}

\paragraph{\textbf{Backward formulation}}

From now on, we give an equivalent backward formulation of the above
stochastic optimal problem (3.2). To do so we need an additional assumption:

(H4) there exists $\alpha>0$ such that $|\sigma(y,u_{1},t)-\sigma
(y,u_{2},t)|\geq\alpha|u_{1}-u_{2}|$ for all $y\in\mathbb{R}^{k},t\in
\lbrack0,T]$ and $u_{1}$, $u_{2}\in\mathbb{R}^{n\times d}.$

Note (H1) and (H4) \ imply the mapping  
$
u\rightarrow\sigma(y,u,t)
$
is a bijection from $\mathbb{R}^{n\times d}$ on to itself for any $(y,t)$.

Let $q\equiv\sigma(y,u,t)$ \ and denote the inverse function by $u=\tilde
{\sigma}(y,q,t)$. Then system (3.1) can be rewritten as%
\[
\left\{
\begin{array}
[c]{rrl}%
-dy(t) & = & f(t,y(t),q(t))dt-q(t)dW_t,\\
y(0) & = & b,\\
-dx(t) & = & F(t,x(t),z(t),y(t),q(t))dt+G(t,x(t),z(t),y(t)
q(t))dW_{t}-z(t)dB_{t},\\
x(0) & = & \xi,
\end{array}
\right.
\]
where $f(t,y,q)=-b(t,y,\tilde{\sigma}(y,q,t))$ and $F(t,x,z
,y,q)=\bar{f}(t,x,z,y,\tilde{\sigma}(y,q,t))$, $G(t,x
,z,y,q)=\bar{g}(t,x,z,y,$ $\tilde{\sigma}(y,q,t))$.

A key observation that inspires our approach of solving problem (3.2) is that,
since $u\rightarrow\sigma(x,u,t)$ is a bijection, $q(\cdot)$ can be regarded
as the control; moreover, by the BSDE theory selecting $q(\cdot)$ is
equivalent to selecting the terminal value $y(T)$. Hence we introduce the
following ``controlled" system:
\begin{equation}
\left\{
\begin{array}
[c]{rrl}%
-dy(t) & = & f(t,y(t),q(t))dt-q(t)dW_t,\\
y(t) & = & \eta,\\
-dx(t) & = & F(t,x(t),z(t),y(t),q(t))dt+G(t,x(t),z(t),y(t),
q(t))dW_{t}-z(t)dB_{t},\\
x(0) & = & \xi,
\end{array}
\right.  \tag{3.3}%
\end{equation}
where the control variables are the random variables $\xi$ and $\eta$ to be
chosen from the following set%
\[
U=\{(\xi,\eta)|\xi\in K_{1},a.s.\ E|\xi|^{2}< \infty,\ \eta\in K_{2}%
,a.s.\ E|\eta|^{2}< \infty.\}.
\]
For each $(\xi,\eta)\in U,$ consider the following cost
\[
J(\xi,\eta)=E[\int_{0}^{T}l(t,x(t),z(t),y(t),q(t))dt+\chi(\xi
)+\lambda(\eta)+\phi(x(T))],
\]
where $l(t,x,z,y,q)=\bar{l}(t,x,z,y_{,}\tilde{\sigma}(y,q,t))$.

This gives rise to the following auxiliary optimization problem:%
\begin{equation}%
\begin{array}
[c]{cl}
& \inf J(\xi,\eta)\\
\text{subject to} & (\xi,\eta)\in U\text{; }y(0)=b.
\end{array}
\tag{3.4}%
\end{equation}
where $y(0)^{(\xi,\eta)}$ is the solution of (3.3) at time $0$ under $\xi$
and $\eta.$

It is clear that the original problem (3.2) is equivalent to the auxiliary one (3.4).

Hence, hereafter we focus ourselves on solving (3.4). The advantage of doing
this is that, since $\xi$ and $\eta$ now are the control variable, the state
constraint in (3.2) becomes a control constraint in (3.4), whereas it is
well-known in control theory that a control constraint is much easier to deal
with than a state constraint. There is, nonetheless, a cost of doing so is
that the original initial condition $y(0)^{(\xi,\eta)}=b$ now becomes a
constraint, as shown in (3.4).

From now on, we denote the solution of (3.3) by $(x^{(\xi,\eta)}%
(\cdot),y^{(\xi,\eta)}(\cdot),q^{(\xi,\eta)}(\cdot),z^{(\xi,\eta)}(\cdot)) $,
whenever necessary, to show the dependence on $(\xi,\eta).$ We also denote
$x^{(\xi,\eta)}(0)$ and $y^{(\xi,\eta)}(0)$ by $x(0)^{(\xi,\eta)}$ and
$y(0)^{(\xi,\eta)}$ respectively. Finally, it is easy to check that
$f,F,G$\ and $l$ satisfy similar conditions in Assumptions (H1)$-$(H3).

We note that this is an special case of (2.4), so by the same method we have
the following result:

\paragraph{\textbf{Stochastic Maximum Principle.}}

We define%
\begin{equation}
\left\{
\begin{array}
[c]{rrl}%
dn(t) & = & (F_{y}^{\ast}(t)m(t)+G_{y}^{\ast}(t)p(t)+f_{y}^{\ast}%
(t)n(t)+l_{y}^{\ast}(t))dt\\
&  & +(F_{q}^{\ast}(t)m(t)+G_{q}^{\ast}(t)p(t)+f_{q}^{\ast}(t)n(t)+l_{q}%
^{\ast}(t))dW_{t},\\
n(0) & = & h_{2},\\
-dm(t) & = & -(F_{x}^{\ast}(t)m(t)+G_{x}^{\ast}(t)p(t)+l_{x}^{\ast}(t))dt\\
&  & +(F_{z}^{\ast}(t)m(t)+G_{z}^{\ast}(t)p(t)+l_{z}^{\ast}(t))dB_{t}%
-p_{t}dW_{t},\\
m(T) & = & -h_{1}\phi_{x}(x^{\ast}(T)),
\end{array}
\right.  \tag{3.5}%
\end{equation}
where $H_{k}^{\ast}(t)=H_{k}(t,x^{\ast}(t),z^{\ast}(t),y^{\ast}(t),q^{\ast
}(t))$ ($H=F,G,l$) and $f_{x}^{\ast}(t)=f_{x}(t,x^{\ast}(t),z^{\ast}(t))$,
$f_{q}^{\ast}(t)=f_{q}(t,x^{\ast}(t),z^{\ast}(t))$ , $h_{1}$ and $h_{2}$ are
defined as in (2.10). It is easy to check that there exist unique processes
$n(\cdot)$ and $(m(\cdot),p(\cdot))$ which solve the above equations.

\begin{theorem}
We assume (H1)$\sim$(H4). Let $(\xi^{\ast},\eta^{\ast})$ be optimal to (3.4)
and $(x^{\ast}(\cdot),y^{\ast}(\cdot),z^{\ast}(\cdot),q^{ \ast}(\cdot))$ be
the corresponding optimal trajectory. Then for arbitrary $(\xi,\eta)\in U,$ we
have the following inequalities hold%
\begin{equation}
\left\{
\begin{array}
[c]{lll}%
\langle m(0)-\chi_{x}(\xi^{\ast}),\xi-\xi^{\ast}\rangle\leq0, &  & \\
\langle n(T)+\lambda_{y}(\eta^{\ast}),\eta-\eta^{\ast}\rangle\geq0. &  &
\end{array}
\right.  \tag{3.6}%
\end{equation}

\end{theorem}

\subsection{System composed of a BDSDE with state constraints}

Although this case describes a controlled BDSDE system with state constraints.
But it seems trivial. Thus, we only give a brief illustration. Given $\eta\in
L^{2}(\Omega,%
\mathcal{F}%
_{T},P;%
\mathbb{R}
^{k})$ and $u(\cdot)\in U[0,T]$, consider the following BDSDE.%
\begin{equation}
\left\{
\begin{array}
[c]{rrl}%
-dy(t) & = & f(t,y(t),q(t),u(t))dt+g(t,y(t),q(t),u(t))dB_{t}%
-q(t)dW_{t},\ 0\leq t\leq T,\\
y(T) & = & \eta.
\end{array}
\right.  \tag{3.7}%
\end{equation}
For given $u(\cdot)\in U[0,T]$ and $f,g$ satisfying (H2) and (H3), from the
Theorem 1.1 in \cite{Pardoux-Peng98}, it is easy to check that there exists a
unique solution $(y(\cdot),q(\cdot))$ of (3.7).

Note that $y(0)^{(\eta,u(\cdot))}$ is a $
\mathcal{F}_{0,T}^{B}$-measurable variable, and
\[
E(h(y^{(\eta,u(\cdot))}(0)))=b.
\]
Now we regard $\eta$ and $u(\cdot)$ as the control variables to be chosen from
the following set%
\[
U=\{(\eta,u(\cdot))|\eta\in K_{1},\ a.s.,\ E[|\eta|^{2}]< \infty,\ u(\cdot)\in
U[0,T]\}.
\]
For each $(\eta,u(\cdot))\in U,$ consider the following cost function%
\[
J(\eta,u(\cdot))=E[\int_{0}^{T}l(t,y(t),q(t),u(t))dt+\lambda(\eta
)+\gamma(y(0))],
\]
which gives rise to the following optimization problem%
\begin{equation}%
\begin{array}
[c]{lll}
&  \qquad \inf\limits_{(\eta,u(\cdot))\in U}J(\eta,u(\cdot))\\
&\text{subject to} \ E(h(y^{(\eta,u(\cdot))}(0)))=b.
\end{array}
\tag{3.8}%
\end{equation}

\paragraph{\textbf{Maximum principle}}

Let%
\[
H(y,q,u,n,\delta)\triangleq\langle f(t,y,q,u),n(t)\rangle+\langle
g(t,y,q,u),\delta(t)\rangle+l(t,y,q,u).
\]
Then the adjoint equation is%
\begin{equation}
\left\{
\begin{array}
[c]{rrl}%
dn(t) & = & (f_{y}^{\ast}(t)n(t)+g_{y}^{\ast}(t)\delta(t)+l_{y}^{\ast
}(t))dt-\delta(t)dB_{t}+(f_{q}^{\ast}(t)n(t) +g_{q}^{\ast}(t)\delta(t)+l_{q}^{\ast}(t))dW_{t},\\
n(0) & = & h_{2}h_{y}(y^{\ast}(0))+h_{0}\gamma_{y}(y^{\ast}(0)),
\end{array}
\right.  \tag{3.9}%
\end{equation}
where $H_{k}^{\ast}(t)=H_{k}(t,y^{\ast}(t),q^{\ast}(t))$ for $H=f,g$; $k=y,q$, respectively.

\begin{theorem}
We assume (H1)$\sim$(H3). Let $\eta^{\ast}$ and $u^{\ast}(\cdot)$ be the
optimal controls and $(y^{\ast}(\cdot),z^{\ast}(\cdot))$ be the corresponding
optimal trajectory. Then we have%
\begin{equation}
\left\{
\begin{array}
[c]{lll}%
\langle n(T)+\lambda_{y}(\eta^{\ast}),\eta-\eta^{\ast}\rangle\geq0, &  & \\
\langle H_{u}(t,y^{\ast},q^{\ast},u^{\ast},n,\delta),u(t)-u^{\ast}%
(t)\rangle\geq0. &  &
\end{array}
\right.  \tag{3.10}%
\end{equation}

\end{theorem}

\subsection{Backward doubly stochastic LQ problem without state constraints}

Consider the following linear system:%
\begin{equation}
\left\{
\begin{array}
[c]{rrl}%
-dx(t) & = & [A(t)x(t)+B(t)z(t)+C(t)y(t)+D(t)q(t)+E(t)u(t)]dt\\
&  & +[A^{\prime}(t)x(t)+B^{\prime}(t)z(t)+C^{\prime}(t)y(t)+D^{\prime
}(t)q(t)+E^{\prime}(t)u(t)]dW_{t} -z(t)dB_{t},\\
x(0) & = & x,\\
-dy(t) & = & [A^{\prime\prime}(t)x(t)+B^{\prime\prime}(t)z(t)+C^{\prime\prime
}(t)y(t)+D^{\prime\prime}(t)q(t)+E^{\prime\prime}(t)u(t)]dt\\
&  & +[A^{\prime\prime\prime}(t)x(t)+B^{\prime\prime\prime}(t)z(t)+C^{\prime
\prime\prime}(t)y(t)+D^{\prime\prime\prime}(t)q(t)+E^{\prime\prime\prime
}(t)u(t)]dB_{t} -q(t)dW_{t},\\
y(T) & = & y,
\end{array}
\right.  \tag{3.11}%
\end{equation}
where $x$ and $y$ are given constants and $H,$ $H^{\prime},$ $H^{\prime\prime
},$ $H^{\prime\prime\prime}$ ($H=A,$ $B,$ $C,$ $D$ and $E$) are corresponding matrixes.

The cost function (2.3) becomes
\[%
\begin{array}
[c]{lll}%
l(t,x(t),z(t),y(t),q(t),u(t)) = & \frac{1}{2}F(t)x(t)\cdot x(t)+\frac{1}%
{2}G(t)z(t)\cdot z(t) & \\
& +\frac{1}{2}H(t)y(t)\cdot y(t)+\frac{1}{2}I(t)q(t)\cdot q(t)+\frac{1}%
{2}J(t)u(t)\cdot u(t), &
\end{array}
\]
and
\[%
\begin{array}
[c]{lll}
&  & \chi(x) = 0,\ \lambda(y) = 0,\ \phi(x) = \frac{1}{2}U(t)x^{2},\ \gamma(y)
= \frac{1}{2}Q(t)y^{2},
\end{array}
\]
where all function of $t$ are bounded and $F(t),G(t),H(t),I(t)$ are symmetric
non-negative define, $J(t),U(t)$ $,Q(t)$ are symmetric uniformly positive definite.

Then from (2.12), the adjoint equations become%
\[
\left\{
\begin{array}
[c]{llll}%
-dm(t) &=& -[A(t)m(t)+A^{\prime}(t)p(t)+A^{\prime\prime}(t)n(t)+A^{\prime
\prime\prime}(t)\delta(t)+F(t)x^{\ast}(t)]dt \\
&&+[B(t)m(t)+B^{\prime}(t)p(t)+B^{\prime\prime}(t)n(t)+B^{\prime
\prime\prime}(t)\delta(t)+G(t)z^{\ast}(t)]dB_{t}-p(t)dW_{t}, \\
m(T) &=& -U(T)x^{\ast}(T),  \\
dn(t)& = &[C(t)m(t)+C^{\prime}(t)p(t)+C^{\prime\prime}(t)n(t)+C^{\prime
\prime\prime}(t)\delta(t)+H(t)y^{\ast}(t)]dt \\
&&+[D(t)m(t)+D^{\prime}(t)p(t)+D^{\prime\prime}(t)n(t)+D^{\prime
\prime\prime}(t)\delta(t)+I(t)q^{\ast}(t)]dW_{t}-\delta(t)dB_{t}, \\
n(0)&=& Q(0)y^{\ast}(0).
\end{array}
\right.
\]
Define%
\[%
\begin{array}
[c]{rl}
& H(t,x,z,y,q,u,m,p,n,\delta)
\triangleq  \langle F(t,x,z,y,q,u),m(t)\rangle+\langle
G(t,x,z,y,q,u),p(t)\rangle\\
& +\langle f(t,x,z,y,q,u),n(t)\rangle+\langle g(t,x,z,y,q,u),\delta
(t)\rangle+l(t,x,z,y,q,u).
\end{array}
\]
Suppose that $K$ is an open set. Then we have the following result%
\[%
\begin{array}
[c]{rl}
& H_{u}(t,x^{\ast},z^{\ast},y^{\ast},q^{\ast},u^{\ast},m,p,n,\delta) =E(t)m(t)+E^{\prime}(t)p(t)+E^{\prime\prime}(t)n(t)+E^{\prime\prime\prime
}(t)\delta(t)+J(t)u(t)=0.
\end{array}
\]
Thus,
\[
u^{\ast}(t)=J^{-1}(t)[E(t)m(t)+E^{\prime}(t)p(t)+E^{\prime\prime
}(t)n(t)+E^{\prime\prime\prime}(t)\delta(t)].
\]
However, the maximum principle gives only the necessary condition for optimal
control. We also have:

\begin{theorem}
The funtion $u^{\ast}(t)=J^{-1}(t)[E(t)m(t)+E^{\prime}(t)p(t)+E^{\prime\prime
}(t)n(t)+E^{\prime\prime\prime}(t)\delta(t)]$ is the unique optimal control
for backward doubly stochastic LQ problems, where $(x^{\ast}(\cdot),z^{\ast
}(\cdot),y^{\ast}(\cdot),q^{ \ast}(\cdot))$ and $(m(\cdot),p(\cdot
),n(\cdot),\delta(\cdot))$ are solutions of above equations.
\end{theorem}

\begin{proof}
First let us prove the $u^{\ast}(\cdot)$ is the optimal control. $\forall
v(\cdot)\in K,$ let $(x^{v}(\cdot),z^{v}(\cdot),$ $y^{v}(\cdot),q^{v}(\cdot))$
be the corresponding trajectory of (3.11), we get%
\[%
\begin{array}
[c]{lll}
& J(v(\cdot))-J(u^{\ast}(\cdot)) & \\
= & \frac{1}{2}\{{\int\nolimits_{0}^{T}} [\langle F(t)x^{v}(t),x^{v}%
(t)\rangle-\langle F(t)x^{\ast}(t),x^{\ast}(t)\rangle+\langle G(t)z^{v}%
(t),z^{v}(t)\rangle & \\
& -\langle G(t)z^{\ast}(t),z^{\ast}(t)\rangle+\langle H(t)y^{v}(t),y^{v}%
(t)\rangle-\langle H(t)y^{\ast}(t),y^{\ast}(t)\rangle+\langle I(t)q^{v}
(t),q^{v}(t)\rangle & \\
& -\langle I(t)q^{\ast}(t),q^{\ast}(t)\rangle+\langle J(t)v(t),v(t)\rangle
-\langle J(t)u^{\ast}(t),u^{\ast}(t)\rangle]dt+\langle U(T)x^{v}%
(T),x^{v}(T)\rangle & \\
& -\langle U(T)x^{\ast}(T),x^{\ast}(T)\rangle+\langle Q(0)y^{v}(0),y^{v}%
(0)\rangle-\langle Q(0)y^{\ast}(0),y^{\ast}(0)\rangle\} & \\
\geq & E\{{\int\nolimits_{0}^{T}} [\langle F(t)x^{\ast}(t),x^{v}(t)-x^{\ast
}(t)\rangle+\langle G(t)z^{\ast}(t),z^{v}(t)-z^{\ast}(t)\rangle & \\
& +\langle H(t)y^{\ast}(t),y^{v}(t)-y^{\ast}(t)\rangle+\langle I(t)q^{\ast
}(t),q^{v}(t)-q^{\ast}(t)\rangle+\langle J(t)u^{\ast}(t),v(t)-u^{\ast
}(t)\rangle]dt & \\
& +\langle U(T)x^{\ast}(T),x^{v}(T)-x^{\ast}(T)\rangle+\langle Q(0)y^{\ast
}(0),y^{v}(0)-y^{\ast}(0)\rangle\}. &
\end{array}
\]
Using Lemma 2.1 to $\langle x^{v}(t)-x^{\ast}(t),m(t)\rangle+\,\langle
y^{v}(t)-y^{\ast}(t),n(t)\rangle$ , we obtain%
\[%
\begin{array}
[c]{rl}
& \langle U(T)x^{\ast}(T),x^{v}(T)-x^{\ast}(T)\rangle+\langle Q(0)y^{\ast
}(0),y^{v}(0)-y^{\ast}(0)\rangle\\
= & -E{\int\nolimits_{0}^{T}} [\langle F(t)x^{\ast}(t),x^{v}(t)-x^{\ast
}(t)\rangle+\langle G(t)z^{\ast}(t),z^{v}(t)-z^{\ast}(t)\rangle\\
& +\langle H(t)y^{\ast}(t),y^{v}(t)-y^{\ast}(t)\rangle+\langle I(t)q^{\ast
}(t),q^{v}(t)-q^{\ast}(t)\rangle\langle J(t)u^{\ast}(t),v(t)-u^{\ast
}(t)\rangle\\
& +\langle E(t)m(t),v(t)-u^{\ast}(t)\rangle+\langle E^{\prime}%
(t)p(t),v(t)-u^{\ast}(t)\rangle+\langle E^{\prime\prime}(t)n(t),v(t)-u^{\ast
}(t)\rangle\\
& +\langle E^{\prime\prime\prime}(t)\delta(t),v(t)-u^{\ast}(t)\rangle]dt.
\end{array}
\]
So by the definition of $u^{\ast}(t),$%
\[%
\begin{array}
[c]{rl}
& J(v(\cdot))-J(u^{\ast}(\cdot))\\
\geq & \langle E(t)m(t),v(t)-u^{\ast}(t)\rangle+\langle E^{\prime
}(t)p(t),v(t)-u^{\ast}(t)\rangle+\langle E^{\prime\prime}(t)n(t),v(t)-u^{\ast
}(t)\rangle\\
& +\langle E^{\prime\prime\prime}(t)\delta(t),v(t)-u^{\ast}(t)\rangle+\langle
J(t)u^{\ast}(t),v(t)-u^{\ast}(t)\rangle]dt
=  0.
\end{array}
\]
From the arbitrariness of $v(\cdot)\in K$, we deduce $u^{\ast}(t)$ is the
optimal control.

The proof of the uniqueness of the optimal control is classical. Assume that
$u^{1}(\cdot)$ and $u^{2}(\cdot)$ are both optimal controls, and the
corresponding trajectories are $(x^{1}(\cdot),z^{ 1}(\cdot),y^{1}(\cdot),q^{
1}(\cdot))$ and $(x^{2}(\cdot),z^{ 2}(\cdot),y^{2}(\cdot),q^{ 2}(\cdot))$. By
the uniqueness of solutions of (3.11), we know that the trajectory
corresponding to $\frac{u^{1}(\cdot)+u^{2}(\cdot)}{2}$ is $(\frac{x^{1}%
(\cdot)+x^{2}(\cdot)}{2},$ $\frac{z^{\ 1}(\cdot)+z^{2}(\cdot)}{2},\frac
{y^{1}(\cdot)+y^{2}(\cdot)}{2},\frac{q^{\ 1}(\cdot)+q^{2}(\cdot)}{2})$, and
notice that $J(t),U(t),Q(t)$ are positive, $F(t),$ $G(t),H(t),I(t)$ are
non-negative, we have
\[
J(u^{1}(\cdot))=J(u^{2}(\cdot))=\alpha\geq0
\]
and
\[%
\begin{array}
[c]{rrl}%
2\alpha & = & J(u^{1}(\cdot))+J(u^{2}(\cdot))\\
& \geq & 2J(\frac{u^{1}(\cdot)+u^{2}(\cdot)}{2})+2E{\displaystyle\int%
\nolimits_{0}^{T}} \langle J(t)\frac{u^{1}(\cdot)-u^{2}(\cdot)}{2},\frac
{u^{1}(\cdot)-u^{2}(\cdot)} {2}\rangle dt\\
& \geq & 2\alpha+\frac{\beta}{2}E{\displaystyle\int\nolimits_{0}^{T}}
|u^{1}(\cdot)-u^{2}(\cdot)|^{2}dt,
\end{array}
\]
here $\beta>0.$ So
\[
E{\displaystyle\int\nolimits_{0}^{T}} |u^{1}(\cdot)-u^{2}(\cdot)|^{2}dt\leq0,
\]
which shows that $u^{1}(\cdot)=u^{2}(\cdot).$
\end{proof}

\textbf{Example:} Consider the following backward doubly stochastic LQ
problem, where $t\in[0,1],$ $u(\cdot)\in U[0,1]=[-1,1]$ and $n=k=d=l=1$.%
\begin{equation}
\left\{
\begin{array}
[c]{rrl}%
-dx(t) & = & \frac{1}{2}(y(t)+2u(t))dW_{t}-z(t)dB_{t},\\
x(0) & = & 0,\\
-dy(t) & = & \frac{1}{2}(x(t)+u(t))dB_{t}-q(t)dW_{t},\\
y(1) & = & 0.
\end{array}
\right.  \tag{3.12}%
\end{equation}
We want to minimize the following cost function
\[
J(u(\cdot))=E\int_{0}^{1}(x(s)^2-y(s)^2+z(s)^2-q(s)^2+2x(s)u(s)-4y(s)u(s))ds+Ex(1)^{2}+Ey(0)^{2}.
\]
From (3.12), we get, for $t\in[0,1],$
\[%
\begin{array}
[c]{lll}%
x(t)=-\int_{0}^{t}(y(s)+2u(s))dW_{s}+\int_{0}^{t}z(s)dB_{s}, &  & \\
y(t)=\int_{t}^{1}(x(s)+u(s))dB_{s}-\int_{t}^{1}q(s)dW_{s}. &  &
\end{array}
\]
By substituting $x(\cdot)$ and $y(\cdot)$ into the cost function, we obtain%
\[
J(u(\cdot))=E\int_{0}^{1}3u(s)^{2}ds.
\]
Thus, the optimal control is $u^{\ast}(t)\equiv0,\ t\in[0,1]$ with the optimal
state trajectory
\[
(x^{\ast}(t),z^{\ast}(t),y^{\ast}(t),q^{\ast}(t))\equiv0,\ t\in[0,1].
\]
The adjoint equations are%
\begin{equation}
\left\{
\begin{array}
[c]{rcl}%
dm(t) & = & n(t) dt+p(t) dW_{t},\\
m(1) & = & 0,\\
dn(t) & = & p(t) dt-\delta(t) dB_{t},\\
n(0) & = & 0.
\end{array}
\right.  \tag{3.13}%
\end{equation}
It is obvious that $(m(\cdot),p(\cdot),n(\cdot),\delta(\cdot))=(0,0,0,0)$ is
the unique solution of the above equation.

\section*{Acknowledgments}

The authors would like to thank Prof. Shige Peng for some useful conversations.

\end{document}